\documentclass[11pt]{article}
\usepackage{graphicx}
\usepackage{hyperref}
\usepackage{color}
\usepackage{amssymb,latexsym,amsmath}
\usepackage{epstopdf}
\newtheorem{lemma}{Lemma}[section]
\newtheorem{theo}[lemma]{Theorem}
\newtheorem{prop}[lemma]{Proposition}

\newtheorem{reslt}[lemma]{Result}
\newtheorem{coro}[lemma]{Corollary}

\newcommand{\proof}{\noindent{\em Proof: }}

\newcommand{\forme}[1]{}

\def\wbull{\hfill\vrule height .9ex width .8ex depth -.1ex}

\begin{document}

\date{\today}
\title{On a characterization of the Grassmann graphs}
\author{{\bf Alexander L. Gavrilyuk}\\
Center for Math Research and Education,
Pusan National University,\\
2, Busandaehak-ro 63beon-gil, Geumjeong-gu, Busan, 46241, Republic of Korea\\
and\\
Krasovskii Institute of Mathematics and Mechanics,\\ 
Kovalevskaya str., 16, Yekaterinburg 620990, Russia\\
e-mail: alexander.gavriliouk@gmail.com\\
\\
{\bf Jack H. Koolen}\\Wen-tsun Wu Key Laboratory of CAS, School of Mathematical Sciences\\ 
University of Science and Technology of China, Hefei 230026, Anhui, PR China\\
e-mail: koolen@ustc.edu.cn}
\maketitle

\begin{abstract}
In 1995, Metsch showed that the Grassmann graph $J_q(n,D)$ of diameter $D\geq 3$ 
is characterized by its intersection numbers with the following possible exceptions: 
\begin{itemize}
\item $n=2D$ or $n=2D+1$, $q\geq 2$; 
\item $n=2D+2$ and $q\in \{2,3\}$;
\item $n=2D+3$ and $q=2$.
\end{itemize}

In 2005, Van Dam and Koolen constructed the twisted Grassmann graphs with the same intersection numbers 
as the Grassmann graphs $J_q(2D+1,D)$, for any prime power $q$ and diameter $D\geq 2$, 
but they are not isomorphic.

We show that the Grassmann graph $J_q(2D,D)$ is characterized 
by its intersection numbers provided that the diameter $D$ is large enough.
\end{abstract}

\section{Introduction}

\noindent {\bf Results}

A finite connected graph $\Gamma$ with vertex set $V(\Gamma)$ and path-length distance 
function $\partial$ is said to be {\it distance-regular} if, for any vertices $x,y\in V(\Gamma)$ 
and any non-negative integers $i,j$, the number $p_{ij}^h$ of vertices at distance $i$ from $x$ 
and distance $j$ from $y$ depends only on $i,j$ and $h:=\partial(x,y)$, and does not depend 
on the particular choice of $x$ and $y$. The numbers $p_{ij}^h$ are called the {\it intersection 
numbers} of $\Gamma$.

Let $\mathbb{F}_q$ be the finite field with $q$ elements, 
and $V$ be the vector space of dimension $n\geq 2$ over $\mathbb{F}_q$.
For an integer $D$, $0<D<n$, 
let $\mathcal{G}_D$ denote the set of all $D$-dimensional subspaces of $V$ 
(i.e., the Grassmannian of $V$). The {\it Grassmann} graph\footnote{The 
Grassmann graphs can be seen as the $q$-``analogues'' of the Johnson graphs, 
which are usually denoted by $J(n,k)$. This explains the notation $J_q(n,D)$.} 
$J_q(n,D)$ 
has $\mathcal{G}_D$ as the vertex set with two vertices being adjacent 
if and only if they meet in a subspace of dimension $D-1$.
As the graphs $J_q(n,D)$ and $J_q(n,n-D)$ are isomorphic (an isomorphism defined 
by mapping each subspace to its orthogonal complement), without loss of generality, 
we further assume that $n\geq 2D$. The Grassmann graph $J_q(n,D)$ 
is distance-regular, and all its intersection numbers are expressed in 
terms of $n,D,$ and $q$.

For a natural number $q\geq 2$, define a function $\chi(q)$ by:
\begin{align}\label{eq-chi}
\chi(q) &= 
\left\{
\begin{matrix}
9 &\mbox{if $q=2,$} \\
8 &\mbox{if $q=3,$} \\
7 &\mbox{if $q\in \{4,5,6\},$} \\
6 &\mbox{if $q\geq 7.$} 
\end{matrix}
\right.
\end{align}

The main result of this paper is as follows.
 
\begin{theo}\label{theo-main}
For a prime power $q$ and a natural number $D\geq \chi(q)$, 
suppose that $\Gamma$ is a distance-regular graph with the same intersection numbers as 
the Grassmann graph  $J_q(2D,D)$. 
Then $\Gamma$ is isomorphic to $J_q(2D,D)$.
\end{theo}

\noindent {\bf Motivation}

A distance-regular graph gives rise to a $P$-polynomial (also known as {\it metric}) 
association scheme and vice versa; it can naturally be seen as 
a finite-analogue of compact $2$-point homogeneous spaces in Riemannian geometry.

It was realized by Delsarte in his seminal work \cite{Delsarte} that 
$P$-polynomial association schemes provide an algebraic framework for 
the theory of error-correcting codes. He introduced 
$Q$-polynomial (also known as {\it cometric}) association schemes
as the dual concept of $P$-polynomial association schemes
(by using the fact that the Bose-Mesner matrix algebra of an association scheme is closed with respect
to both the standard and entry-wise products, where the latter one is in 
a sense ``dual'' to the former one), 
and showed that they provide an algebraic framework for the theory of combinatorial designs.
This unified coding theory and design theory, for which $P$- and $Q$-polynomial
association schemes serve as underlying spaces, respectively, and led to the Delsarte theory, 
a linear-algebraic approach to their problems.

It was further observed in the introduction of the monograph by Bannai and Ito 
\cite{BI} that the theory of designs in $Q$-polynomial 
association schemes goes in parallel with that of combinatorial configurations in compact symmetric 
spaces of rank $1$ (for example, spherical designs). From this point of view, $Q$-polynomial 
association schemes can be seen as a finite-analogue of compact symmetric spaces of rank $1$. 

Recall that compact symmetric spaces of rank $1$ were classified by Cartan \cite{Cartan}, 
and it was shown by Wang \cite{Wang} that a compact symmetric space of rank $1$ 
is a compact $2$-point homogeneous space and vice versa. These fundamental results 
from Riemannian geometry lead to the following conjecture proposed in \cite{BI}: 
``Primitive $P$-polynomial association schemes of sufficiently 
large diameter are $Q$-polynomial, and vice versa'' and to the problem of classification 
of primitive $(\text{both~}P$ and $Q)$-polynomial association schemes, which would be a  
finite-analogue of Cartan's classification. 
The list of currently known examples of primitive $(\text{both~}P$ and $Q)$-polynomial 
association schemes includes 20 families of unbounded diameter, and most of them 
arise from classical algebraic objects such as dual polar spaces and forms over finite fields \cite{BI}.
We refer the reader to \cite{BI}, \cite{BBI}, \cite{SurveyDRG} for the current state of the problem 
and the detailed description of the known examples.

The set $\mathcal{G}_D$ gives rise to the {\it Grassmann} association scheme 
whose relations (classes) $R_0,R_1,\ldots,R_D$ are given by:
\begin{equation}\label{eq-grassmanRi}
(U,W)\in R_i \Leftrightarrow \mathsf{dim}(U\cap W)=D-i,~~(U,W\in \mathcal{G}_D).
\end{equation}
This scheme is $(\text{both~}P$ and $Q)$-polynomial and, moreover, 
its relations arise as orbits of a transitive on $\mathcal{G}_D$ group $P\Gamma L(V)$ 
whose action is componentwise extended to $\mathcal{G}_D\times \mathcal{G}_D$. 
In particular, this means that the {\it Grassmann graph} $J_q(n,D)$, which is defined by relation $R_1$, 
is distance-transitive.

The classification problem requires a characterization of the known examples of 
$(\text{both~}P$ and $Q)$-polynomial association schemes (i.e., $Q$-{\it polynomial} 
distance-regular graphs) by their intersection numbers.
Such a characterization was shown 
for the following families of association schemes: 
the Hamming schemes \cite{Egawa}, the Johnson schemes \cite{TerwilligerJohnson},  
and their quotients \cite{Neumaier,Buss,Metsch971,Metsch03,Metsch97,GavrKoolen}, 
the schemes of Hermitian forms \cite{IVHer,KiteFree} and the schemes of dual polar spaces 
of unitary type (in even dimension) \cite{IVA2d,BWil}, the association schemes of bilinear forms 
\cite{CuypersBF,Huang,Metsch95,Bil2d} (some cases left open).
This paper contributes to the solution of the classification problem in regards 
to the Grassmann schemes.

\noindent {\bf Previous works}

Much attention has been paid to the problem of characterization of the Grassmann 
graphs (or, in terms of finite geometry, the Grassmann manifolds as a class of incidence 
structures satisfying certain conditions), see 
\cite{Bichara,Biondi,Cohen,CohenCoop,Cooperstein,CuypersJq,Lore,Melone,RCS,Shult,Shult2,Sprague,Tallini}.
The strongest result in this direction was obtained by Metsch in \cite{Metsch95}, 
where he showed that the Grassmann graph $J_q(n,D)$, $D\geq 3$, can be uniquely determined 
as a distance-regular graph by its intersection numbers unless one of the following few cases holds: 
\begin{itemize}
\item $n=2D$ or $n=2D+1$, $q\geq 2$;
\item $n=2D+2$ and $q\in \{2,3\}$;
\item $n=2D+3$ and $q=2$.
\end{itemize}

Note that a characterization of the Grassmann graphs $J_q(n,2)$ in terms of parameters 
is not possible, as these graphs have the same intersection numbers as the block graphs 
of 2-designs with parameters $(v,k,\lambda)=(\frac{q^n-1}{q-1},q+1,1)$, 
and there exist many pairwise non-isomorphic of those \cite{Jung,Wilson}.

The result of Metsch relies on a characterization 
of the incidence structure formed by the vertices and the maximum cliques of 
the Grassmann graph $J_q(n,D)$.
Recall that a {\it partial linear space} is an incidence structure $(P,L,I)$, 
where $P$ and $L$ are sets (whose elements are called {\it points} and {\it lines}, respectively) 
and $I\subseteq P\times L$ is the {\it incidence relation} such that every line is incident 
with at least two points and there exists at most one line through any two distinct points.
The {\it point graph} of the incidence structure $(P,L,I)$ is a graph defined on $P$ as 
the vertex set, with two points being adjacent if they are collinear. 
Observe that $J_q(n,D)$ has two families of maximal cliques corresponding to the sets 
$\mathcal{G}_{D-1}$ and $\mathcal{G}_{D+1}$: 
the maximal cliques of the first family are the collections of $D$-subspaces of $V$ containing 
a fixed subspace of dimension $D-1$, 
and each of them is of size $(q^{n-(D-1)}-1)/(q-1)$, while the maximal cliques of the other family
are the collections of $D$-subspaces of $V$ contained in a fixed subspace of dimension $D+1$, 
and each of them is of size $(q^{D+1}-1)/(q-1)$. 
Every edge of $J_q(n,D)$ is contained in a unique clique of each family, 
and one can then see that $(\mathcal{G}_D,\mathcal{G}_{D-1},\supset)$ is a partial linear space 
with the point graph isomorphic to $J_q(n,D)$.

In \cite{RCS}, Ray-Chaudhuri and Sprague characterized $(\mathcal{G}_D,\mathcal{G}_{D-1},\supset)$ as 
a class of partial linear spaces satisfying certain regularity conditions.
Let $\Gamma$ denote a distance-regular graph with the same intersection numbers as $J_q(n,D)$. 
A key idea of Metsch \cite{Metsch95} was to construct a partial linear space from $\Gamma$ 
by taking its vertices as the points and a set $\mathcal{L}$ of (sufficiently large but 
not necessarily maximum) cliques as the lines, and then, by using the result of Ray-Chaudhuri and Sprague, 
to show that this incidence structure is isomorphic to $(\mathcal{G}_D,\mathcal{G}_{D-1},\supset)$.
Indeed, if 
every edge of $\Gamma$ is contained in a unique clique of $\mathcal{L}$, then 
$(V(\Gamma),\mathcal{L},\in)$ is a partial linear space and $\Gamma$ is its point graph.
To construct large cliques and to show the existence of such a set $\mathcal{L}$, 
Metsch used a counting technique known as a Bose-Laskar type argument \cite{Metsch91}.

This approach fails in the open cases mentioned above (in particular, when $n=2D$, 
the cliques of both families in $J_q(n,D)$ have the same size, and so every edge 
is contained in two maximum cliques). Moreover, Van Dam and Koolen \cite{vDK} discovered 
a new family of distance-regular graphs, the so-called {\it twisted Grassmann} graphs, 
which have the same intersection numbers as $J_q(2D+1,D)$ for any prime power $q$ 
and which are not point graphs of any partial linear space. 
(For the detailed study of these graphs, see \cite{TwistG}, \cite{TwistG2}, 
\cite{TwistGr3}, \cite{TwistGr4}.) 
This demonstrates that that classification problem of the remaining open cases is very challenging.

\noindent {\bf Outline of the proof}

We recall that the intersection numbers of most of known primitive $Q$-polynomial distance-regular graphs 
(in particular, of those related to classical groups and groups of Lie type) 
can be expressed in terms of the so-called {\it classical parameters}, namely, 
the diameter $D$ and three other parameters $b,\alpha$, and $\beta$ 
(see Section \ref{ssect-Qpoly}). 
In this paper we prove the following theorem, which 
implies Theorem \ref{theo-main} immediately.

\begin{theo}\label{theo-main-1}
Suppose that $\Gamma$ is a distance-regular graph  
with classical parameters $(D,b,\alpha,\beta) = (D,q,q,\frac{q^{D+1}-1}{q-1}-1)$ 
for some natural number $q\geq 2$.
If $D\geq \chi(q)$, then $q$ is a prime power and $\Gamma$ is isomorphic to the Grassmann graph 
$J_q(2D,D)$.
\end{theo}

Together with the result of Metsch \cite{Metsch95}, this yields the following.

\begin{coro} 
Suppose that $\Gamma$ is a distance-regular graph  
with classical parameters $(D,b,\alpha,\beta) = (D,q,q,\frac{q^{n-D+1}-1}{q-1}-1)$ 
for some natural numbers $q\geq 2$, $n\geq 2D$, $D\geq 3$.
If one of the following conditions holds:
\begin{itemize}
\item $n=2D$ and $D\geq \chi(q)$,
\item $n\geq 2D+\mathsf{max}(6-q,2)$,
\end{itemize}
then $q$ is a prime power and $\Gamma$ is isomorphic to the Grassmann graph $J_q(n,D)$.
\end{coro}


The proof of Theorem \ref{theo-main-1} exploits several different techniques 
and relies on a characterization of another partial linear space derived 
from the Grassmann graph $J_q(n,D)$ whose points again are the vertices and whose lines 
are the singular lines (here by a {\it singular line} we mean the non-trivial intersection 
of two cliques from different families). Namely, combining a characterization of such partial linear spaces 
obtained by Cooperstein \cite{Cooperstein} and Cohen \cite{Cohen} with a work of Numata \cite{Numata} 
allows to recognize the Grassmann graphs by their {\it local graphs}, i.e., the subgraphs induced 
by the neighbours of vertices. We thus call this result the Numata-Cohen-Cooperstein theorem 
(see \cite[Theorem~9.3.8]{BCN} and Section \ref{sect-local}). 

In order to recover the local structure of a graph $\Gamma$ with the same intersection numbers 
as the Grassmann graph $J_q(2D,D)$, we first use two deep consequences of the Terwilliger algebra 
theory: the {\it triple intersection numbers} (see Section \ref{sect-triple}) and 
the {\it Terwilliger polynomial} (see Section \ref{sect-Tpoly}). 
The former one restricts a possible structure of the local graphs, 
while the latter one restricts possible eigenvalues of their adjacency matrices. 
In Section \ref{sect-prelim}, by using these two ingredients, we show 
that the local graphs of $\Gamma$ share many properties with the local graphs of $J_q(2D,D)$, 
in particular, their adjacency matrices have exactly the same spectrum.

Although in general it appears to be a hard problem to recognize a graph from its spectrum 
\cite[Chapter~14]{BH},
we proceed in Section \ref{sect-bigcliques} by proving that the local graphs of $\Gamma$ 
are indeed isomorphic to those of $J_q(2D,D)$ provided that the diameter $D$ is not too small. 
The proof of this step is based on \cite{Riaz}, \cite{Aida}, and it combines 
some tricks from algebraic graph theory (see Section \ref{sect-spectra}) with a counting argument 
in order to construct large cliques in the local graph (which mimics the proof by Metsch).


Thus, the present work settles the problem of characterization of 
the Grassmann graphs $J_q(n,D)$ in the case $n=2D$ and the diameter $D$ is large enough 
(although our results in Section \ref{sect-prelim} provide certain evidence that 
there should not be exceptions like the twisted Grassmann graphs even when the diameter is small). 
For the cases $n=2D+2$ and $n=2D+3$, a characterization of the Grassmann graphs $J_2(n,D)$ 
will be shown in a forthcoming paper of the first author.

\section{Basic theory}\label{sect-basic}

The main purpose of this section is to recall and to fix some basic terminology 
and notation from algebraic graph theory. 
For more comprehensive background on distance-regular graphs and association
schemes, we refer the reader to \cite{BI}, \cite{BCN}, and \cite{SurveyDRG}.

\subsection{Graphs and their eigenvalues}\label{sect-spectra}

All graphs considered in this paper are finite and simple. 
Let $\Gamma$ be a connected graph. 
The distance $\partial(x,y):=\partial_{\Gamma}(x,y)$ between any two vertices $x,y$ of
$\Gamma$ is the length of a shortest path connecting $x$ and $y$ in $\Gamma$. 
For a subset $X$ of the vertex set of $\Gamma$, we will also 
write $X$ for the subgraph of $\Gamma$ induced by $X$.
For a vertex $x\in \Gamma$, define $\Gamma_i(x)$ 
to be the set of vertices that are at distance precisely $i$ from $x$ ($0\leq i\leq D$),
where $D:=\mathsf{max}\{\partial(x,y)\mid x,y\in \Gamma\}$ is the {\it diameter} of $\Gamma$. 
In addition, define $\Gamma_{-1}(x)=\Gamma_{D+1}(x)=\emptyset$. 
The subgraph induced by $\Gamma_1(x)$ is called the {\it neighborhood} or 
the {\it local graph} of a vertex $x$.
We often use $\Gamma(x)$ instead of $\Gamma_1(x)$ for short, and we write $x\sim_{\Gamma} y$ 
or simply $x\sim y$ if two vertices $x$ and $y$ are adjacent in $\Gamma$. 
A graph $\Gamma$ is {\it regular} with {\it valency} $k$ if 
the local graph $\Gamma(x)$ contains precisely $k$ vertices for all $x\in \Gamma$.

For a set $\{x_1,x_2,\ldots,x_s\}$ of vertices of $\Gamma$, let 
$\Gamma(x_1,x_2,\ldots,x_s)$ denote $\cap_{i=1}^s \Gamma(x_i)$.
In particular, for a pair $x,y$ of vertices of $\Gamma$ with $\partial(x,y)=2$, 
the subgraph induced on $\Gamma(x,y)$ is commonly known as the $\mu$-{\it graph} (of $x$ and $y$).

The {\it eigenvalues} of a graph $\Gamma$ are the eigenvalues of its adjacency matrix $A:=A(\Gamma)$.
If, for an eigenvalue $\eta$ of $\Gamma$, its eigenspace contains 
a vector orthogonal to the all-one vector, we say that $\eta$ is {\it non-principal}.
If $\Gamma$ is regular with valency $k$, then all its eigenvalues are non-principal 
unless the graph is connected and then the only eigenvalue that is principal 
is its valency $k$.

Let $\Gamma$ be a graph on $v$ vertices with spectrum 
$[\theta_0]^{m_{0}}$, $[\theta_1]^{m_1}$, $\ldots$, $[\theta_d]^{m_d}$, 
where $\theta_0>\theta_1>\ldots>\theta_d$ are all distinct eigenvalues of $\Gamma$, 
and $m_0,m_1,\ldots,m_d$ are their respective multiplicities. 
Then
\begin{equation*}\label{eq-trace}
\mathsf{Tr}(A^{\ell}) = \sum_{i=0}^{d}m_i \theta_i^{\ell} =  \mbox{the~number~of~closed~walks~of~length~$\ell$~in~$\Gamma$}
~~(\ell \geq 0)
\end{equation*}
where $\mathsf{Tr}(A^{\ell})$ is the trace of matrix $A^{\ell}$ (cf. \cite[Lemma 2.5]{biggs}), 
so that 
\begin{equation}\label{eq-trA-1}
\mathsf{Tr}(A^0)=\sum_{i=0}^{d}m_i=v,~~\mathsf{Tr}(A)=0,
\end{equation}
and, if $\Gamma$ is regular with valency $k$, then $\theta_0=k$, and for $\ell=2$ we obtain:
\begin{equation}\label{eq-trA-2}
\mathsf{Tr}(A^2)=vk. 
\end{equation}

Suppose that a connected graph $\Gamma$ has just $4$ distinct eigenvalues 
and it is regular with valency $k$. Then its adjacency matrix $A$ satisfies \cite{Hoffman,vD95}:
\[
A^3-\big(\sum_{i=1}^3\theta_i\big)A^2+\big(\sum_{1\leq i<j\leq 3}\theta_i\theta_j\big)A
-\theta_1\theta_2\theta_3I=\frac{\prod_{i=1}^3(k-\theta_i)}{v}J,
\]
hereinafter $I$ is the identity matrix, and $J$ is the all-one matrix. 
This shows that $A^3$ has a constant diagonal, 
and thus so does $A^{\ell}$, $\ell=4,5,\ldots$, 
which implies the following result (see \cite{vD95}).

\begin{reslt}\label{reslt-vanDam}
With the above assumption, the following holds.
\begin{itemize}
\item[$(1)$] The number of triangles through any vertex of $\Gamma$ equals
$\frac{1}{2v}\mathsf{Tr}(A^3)$.
\item[$(2)$] The number of quadrangles (a quadrangle may have diagonal edges) 
through any vertex of $\Gamma$ equals
$\frac{1}{2v}\mathsf{Tr}(A^4)-k^2+k/2$.
\end{itemize} 
\end{reslt}
%


Recall that an $s$-\emph{clique} of a graph is its complete subgraph 
(i.e., every two of its vertices are adjacent) with exactly $s$ vertices. 
We call an $s$-clique simply a \emph{clique} if we do not refer to its cardinality.
By the $(s\times t)$-{\it grid}, we mean the Cartesian product of 
two complete graphs on $s$ and $t$ vertices, which is also isomorphic to 
the {\it line} graph of a complete bipartite graph with parts of size $s$ and $t$.
In particular, the $(s\times s)$-grid has spectrum
\begin{equation}\label{eq-gridspectr}
[2(s-1)]^1, [s-2]^{2(s-1)}, [-2]^{(s-1)^2},
\end{equation}
and, moreover, any graph with this spectrum is the $(s\times s)$-grid 
unless $s=4$ (since the Shrikhande graph has the same spectrum as the $(4\times 4)$-grid, 
see \cite{Shrikhande}).
In general, we say that two graphs are {\it cospectral} if they have the same spectrum.

A graph $\Gamma$ is said to be the $q$-{\it clique extension} of a graph $\Delta$ if there 
exists a mapping $\varepsilon$ of the vertex set of $\Gamma$ onto the vertex set of $\Delta$ 
such that $|\varepsilon^{-1}(x)|=q$ for every $x\in \Delta$ and two distinct vertices 
$u,w\in \Gamma$ are adjacent if and only if their images $\varepsilon(u)$ and $\varepsilon(w)$ 
are either equal or adjacent in $\Delta$.
If $A$ is the adjacency matrix of $\Delta$, then the adjacency matrix of $\Gamma$ 
can be written as $J_{q\times |\Delta|}\otimes (A+I_{|\Delta|})-I_{q\times |\Delta|}$
(where $\otimes$ is the Kronecker product),  
whence one can see the following.

\begin{reslt}\label{reslt-cliqueextens}
Suppose that, for an integer $q\geq 1$, a graph $\Gamma$ is the $q$-clique extension of a graph $\Delta$.
Then, for each eigenvalue $\theta$ with $\theta\ne -1$ of $\Delta$, $(q(\theta + 1) - 1)$ is an eigenvalue 
of $\Gamma$ with the same multiplicity. All other eigenvalues of $\Gamma$ are equal to $-1$.
\end{reslt}

We recall one more important result from algebraic graph theory, 
which will be referred to as {\it interlacing}, see (\cite[Section~2.5]{BH}).

\begin{reslt}\label{lem-interlacing}
Let $N$ be a real symmetric $n\times n$ matrix with eigenvalues $\theta_1\geq\ldots\geq\theta_n$. For some
$m<n$, let $R$ be a real $n\times m$ matrix with orthonormal columns, i.e., $R^{\top}R=I$, and 
let $M:=R^{\top}NR$ have eigenvalues $\eta_1\geq\ldots\geq\eta_m$. Then
  the eigenvalues of $M$ interlace those of $N$, i.e., 
  \begin{equation*}
  \theta_i\geq\eta_i\geq\theta_{n-m+i},~~\mbox{for~}i=1,\ldots,m.
  \end{equation*}
\end{reslt}

In particular, this result applies to any principal submatrix $M$ of $N$, as 
one can choose $R$ to be permutation-similar to
$\left(
  \begin{array}{c}
    I \\
    O \\
  \end{array}
\right)$.

Further, let $\pi:=\{V_1,\ldots,V_m\}$ be a partition of the set of columns of a real symmetric matrix $N$ 
and let $N$ be partitioned according to $\pi$ as
\begin{equation*}
\left(
  \begin{array}{ccc}
    N_{1,1} & \ldots & N_{1,m} \\
    \vdots & \ddots & \vdots \\
    N_{m,1} & \ldots & N_{m,m} \\
  \end{array}
\right),
\end{equation*}
where $N_{i,j}$ denotes the submatrix (block) of $N$ formed by columns in $V_j$ 
and by rows that correspond to columns in $V_i$.
The \emph{characteristic matrix} $P$ of $\pi$ is the $n\times m$ matrix whose $j$th column
is the characteristic vector of $V_j~(j=1,\ldots,m)$.
The \emph{quotient matrix} of $N$ with respect to $\pi$ is the $m\times m$ matrix $Q$ whose entries 
are the average row sums of the blocks $N_{i,j}$ of $N$, i.e., 
\begin{equation*}
(Q)_{i,j}=\frac{1}{|V_i|}(P^{\top}NP)_{i,j}.
\end{equation*}


\begin{reslt}\label{lema-matrix-interlacing}
Let $N$ be a real symmetric matrix, and $\pi$ a partition of the set of its columns. 
  Then the eigenvalues of 
  the quotient matrix of $N$ with respect to $\pi$ 
  interlace those of $N$.
\end{reslt}


\subsection{Distance-regular graphs}\label{sect-drg}


A connected graph $\Gamma$ of diameter $D$ is called {\it distance-regular} 
if there exist integers $b_i$ and $c_i$, $0\leq i\leq D$, such that, 
for any pair of vertices $x,y\in \Gamma$ with $\partial(x,y)=i$, there are precisely 
$c_i$ neighbours of $y$ in $\Gamma_{i-1}(x)$ and $b_i$ neighbours of $y$ in $\Gamma_{i+1}(x)$. 
In particular, a distance-regular graph is regular with valency $k:=b_0=c_i+a_i+b_i$. 
We define $a_i:=k-b_i-c_i$, $1\leq i\leq D$, and note that 
$a_i=|\Gamma(y)\cap \Gamma_i(x)|$ holds for any pair of vertices $x,y$ with 
$\partial(x,y)=i$. 
We also define $k_i:=\frac{b_0\cdots b_{i-1}}{c_1\cdots c_i}$, $1\leq i\leq D$, 
and note that $k_i=|\Gamma_i(x)|$ for all $x\in \Gamma$ (so that $k=k_1$).
The array $\{b_0,b_1,\ldots,b_{D-1};c_1,c_2,\ldots,c_D\}$ is 
called the {\it intersection array} of the distance-regular graph $\Gamma$.

A graph $\Gamma$ is distance-regular if and only if, for all integers $h,i,j$ 
with $0\leq h,i,j\leq D$ and all vertices $x,y\in \Gamma$ with $\partial(x,y)=h$, 
the number 
\begin{equation*}\label{eq-pijk}
p^h_{ij}:=|\{z\in \Gamma\mid \partial(x,z)=i,~\partial(y,z)=j\}|=|\Gamma_i(x)\cap \Gamma_j(y)|
\end{equation*}
does not depend on the particular choice of $x,y$. The numbers $p^h_{ij}$ 
are called the {\it intersection numbers} of $\Gamma$.
Note that $k_i=p^0_{ii}$, $c_i=p^i_{1i-1}$, $a_i=p^i_{1i}$, and $b_{i-1}=p^{i-1}_{1i}$, 
$1\leq i\leq D$, and all intersection numbers $p_{ij}^h$ can be calculated 
from the intersection array of $\Gamma$, see \cite[Lemma~4.1.7]{BCN}.


%


\subsection{The Bose-Mesner algebra}\label{ssect-BMalgebra}
Let $\Gamma$ be a distance-regular graph of diameter $D$. 
For each integer $i$ with $0\leq i\leq D$, define the $i$th {\it distance matrix} $A_i$ 
of $\Gamma$ whose rows and columns are indexed by the vertex set of $\Gamma$, and, 
for any $x,y\in \Gamma$, 
\begin{equation*}
(A_i)_{x,y} = \left \{ \begin{aligned}
1\text{~if~}\partial(x,y)=i,\\
0\text{~if~}\partial(x,y)\ne i.
\end{aligned}\right.
\end{equation*}

Then $A:=A_1$ is the {\it adjacency matrix} of $\Gamma$, 
$A_0=I$, $A_i^{\top}=A_i$ ($0\leq i\leq D$), and 
\begin{equation*}
A_iA_j=\sum_{h=0}^D p^h_{ij}A_h ~~~ (0\leq i,j\leq D),
\end{equation*}
in particular,
\begin{equation*}
AA_i=b_{i-1}A_{i-1}+a_{i}A_{i}+c_{i+1}A_{i+1} ~~~~ (1\leq i\leq D-1),
\end{equation*}
\begin{equation*}
AA_D=b_{D-1}A_{D-1}+a_{D}A_{D}, 
\end{equation*}
and this implies that $A_i=p_i(A)$ for certain polynomial $p_i$ of degree $i$ for $0\leq i\leq D$. 

The {\it Bose-Mesner} algebra $\mathcal{ M}$ of $\Gamma$ is the matrix algebra generated 
by $A$ over ${\mathbb R}$. It follows that $\mathcal{ M}$ has dimension $D+1$, 
and it is spanned by the set of matrices $A_0=I,A_1,\ldots,A_D$, which form a basis of $\mathcal{ M}$.
Since the algebra $\mathcal{ M}$ is semi-simple and commutative, $\mathcal{ M}$ also has 
a basis of pairwise orthogonal idempotents $E_0:=\frac{1}{|\Gamma|}J,E_1,\ldots,E_D$ 
(the so-called {\it primitive idempotents} of $\mathcal{ M}$) satisfying:
\begin{equation*}
E_iE_j=\delta_{ij}E_i~~(0\leq i,j\leq D),~~E_i=E_i^{\top}~~(0\leq i\leq D), 
\end{equation*}
\begin{equation*}
E_0+E_1+\cdots+E_D=I.
\end{equation*}
 
A distance-regular graph of diameter $D$ has precisely $D+1$ distinct 
eigenvalues, which can be calculated from its intersection array, see 
\cite[Section 4.1.B]{BCN}. In fact, $E_j$ ($0\leq j\leq D$) turns out to be the matrix of rank 
$m_j:=\mathsf{Tr}(E_j)$
representing orthogonal projection onto the eigenspace of $A$ corresponding to some eigenvalue, 
say $\theta_j$, with multiplicity $m_j$ of $\Gamma$. In other words, one can write
\begin{equation*}
A=\sum_{j=0}^D \theta_jE_j,
\end{equation*} 
where $\theta_j$ ($0\leq j\leq D$) are the real and pairwise distinct 
scalars, which are exactly the eigenvalues of $\Gamma$.
We say that the eigenvalues (and the corresponding idempotents $E_0,E_1,\ldots,E_D$) 
are in {\it natural} order if $b_0=\theta_0>\theta_1>\ldots>\theta_D$.

%

The Bose-Mesner algebra $\mathcal{ M}$ is also closed under entrywise 
matrix multiplication, denoted by $\circ$. The matrices $A_0$, $A_1$, $\ldots$, $A_D$ 
are the primitive idempotents of $\mathcal{ M}$ with respect to $\circ$, i.e., 
$A_i\circ A_j= \delta_{ij}A_i$, and $\sum_{i=0}^D A_i=J$. 
This implies that
\[
E_i\circ E_j=\sum_{h=0}^{D} q_{ij}^h E_h ~~~ (0\leq i,j\leq D)
\]
holds for some real numbers $q_{ij}^h$, known as the {\it Krein parameters} of $\Gamma$. 

The Krein parameters $q_{ij}^h$ can be seen as a counterpart to the intersection numbers $p_{ij}^h$, 
however, they do not have to be integers and do not satisfy, in general, the triangle inequality  
as $p_{ij}^h$ do (i.e., $p_{ij}^h=0$ whenever $i+j<h$ or $|i-j|>h$).

\subsection{$Q$-polynomial distance-regular graphs and classical parameters}\label{ssect-Qpoly}
Let $\Gamma$ be a distance-regular graph of diameter $D$, and $E$ be one of the primitive idempotents 
of its Bose-Mesner algebra. The graph $\Gamma$ is called {\it $Q$-polynomial} with respect to $E$ 
(or with respect to an eigenvalue $\theta$ of $A$ corresponding to $E$) 
if there exist real numbers $c_i^*$, $a_i^*$, $b_{i-1}^*$ ($1\leq i\leq D$) 
and an ordering of the primitive idempotents $E_0,E_1,\ldots,E_{D}$
such that $E_0=\frac{1}{|\Gamma|}J$ and $E_1=E$, and
\[
E_1\circ E_i=b_{i-1}^*E_{i-1} + a_i^*E_i + c_{i+1}^*E_{i+1} ~~~ (1\leq i\leq D-1),
\]
\[
E_1\circ E_D=b_{D-1}^*E_{D-1} + a_D^*E_D.
\]

We call such an ordering of primitive idempotents (and that of the corresponding eigenvalues of $\Gamma$) 
$Q$-{\it polynomial}. Note that a $Q$-polynomial ordering of the eigenvalues/idempotents does not 
have to be the natural one.
One can see that, in terms of the Krein parameters, 
$m_i=q^0_{ii}$, $c_i^*=q^i_{1i-1}$, $a_i^*=q^i_{1i}$, 
and $b_{i-1}^*=q^{i-1}_{1i}$ for $1\leq i\leq D$ 
(here we use index $i$ with respect to a $Q$-polynomial ordering of the primitive idempotents).
In this case, the Krein parameters satisfy the triangle inequality, i.e.,  
$q_{ij}^h=0$ whenever $i+j<h$ or $|i-j|>h$.



Recall that the $q$-{\it ary Gaussian binomial coefficient} is defined by
\[
\left[
\begin{matrix} 
n \\ 
m
\end{matrix}
\right]_q=
\frac{(q^n-1)(q^{n-1}-1)\cdots (q^{n-m+1}-1)}{(q^m-1)(q^{m-1}-1)\cdots (q-1)}.
\]

We say that a distance-regular graph $\Gamma$ of diameter $D$ has {\it classical parameters} 
$(D,b,\alpha,\beta)$ 
if the intersection numbers of $\Gamma$ 
satisfy 
\begin{equation}\label{classparamc_i}
c_i=\genfrac{[}{]}{0pt}{}{i}{1}\Big(1+\alpha\genfrac{[}{]}{0pt}{}{i-1}{1}\Big),
\end{equation}
\begin{equation}\label{classparamb_i}
b_i=\Big(\genfrac{[}{]}{0pt}{}{D}{1}-\genfrac{[}{]}{0pt}{}{i}{1}\Big)\Big(\beta-\alpha\genfrac{[}{]}{0pt}{}{i}{1}\Big),
\end{equation}
where 
\begin{equation}\label{eq-classical-qbinom}
\genfrac{[}{]}{0pt}{}{j}{1}:=\genfrac{[}{]}{0pt}{}{j}{1}_b=1+b+b^2+\cdots+b^{j-1}.
\end{equation}

Note that a distance-regular graph with classical parameters is $Q$-polynomial, 
see \cite[Corollary 8.4.2]{BCN}. 
By \cite[Table~6.1,~Theorem~9.3.3]{BCN}, we have the following result.

\begin{reslt}\label{reslt-classparam}
The Grassmann graph $J_q(n,D)$, $n\geq 2D$, has classical parameters 
\[
(D,b,\alpha,\beta) = (D,q,q,\genfrac{[}{]}{0pt}{}{n-D+1}{1}_q-1).
\]

A distance-regular graph with these classical parameters has 
intersection array given by ($1\leq j \leq D$)
\begin{equation}\label{eq-bj}
b_{j-1}=q^{2j-1}
\left[
\begin{matrix} 
n-D-j+1 \\ 
1
\end{matrix}
\right]_q
\left[
\begin{matrix} 
D-j+1 \\ 
1
\end{matrix}
\right]_q,
\end{equation}
\begin{equation}\label{eq-cj}
c_j=
\left[
\begin{matrix} 
j \\ 
1
\end{matrix}
\right]^2_q,
\end{equation}
and its eigenvalues and their respective multiplicities are given by
(for $0\leq j \leq D$)
\begin{equation}\label{eq-thetaj}
\theta_{j}=q^{j+1}
\left[
\begin{matrix} 
n-D-j \\ 
1
\end{matrix}
\right]_q
\left[
\begin{matrix} 
D-j \\ 
1
\end{matrix}
\right]_q
-
\left[
\begin{matrix} 
j \\ 
1
\end{matrix}
\right]_q,
\end{equation}
\begin{equation}\label{eq-mj}
m_j=
\left[
\begin{matrix} 
n \\ 
j
\end{matrix}
\right]_q
-
\left[
\begin{matrix} 
n \\ 
j-1
\end{matrix}
\right]_q.
\end{equation}
\end{reslt}

\section{On the Terwilliger algebra of a $Q$-polynomial distance-regular graph}

In this section, we explain two key ingredients of our proof, which are based 
on the Terwilliger algebra theory: triple intersection numbers and the Terwilliger polynomial.

\subsection{Triple intersection numbers}\label{sect-triple}

Let $\Gamma$ denote a distance-regular graph of diameter $D\geq 3$. 
Pick any 3-tuple $xyz$ of vertices of $\Gamma$ such that $y$ and $z$ 
are neighbours of $x$. Let $[\ell,m,n]:=[\ell,m,n]_{x,y,z}$ denote 
the {\it triple intersection number} (with respect to $xyz$) defined by:
\[
[\ell,m,n]:=|\Gamma_{\ell}(x)\cap \Gamma_{m}(y)\cap \Gamma_{n}(z)|.
\]

Unlike the intersection numbers, 
the triple intersection numbers $[\ell,m,n]$ depend, in general, 
on the choice of $x,y,z$. 
On the other hand, it is known that vanishing of some of the Krein parameters 
of a distance-regular graph
often leads to non-trivial equations involving triple intersection numbers as the unknowns, 
see, for example, \cite[Theorem~2.3.2]{BCN}, \cite{JurisicCoolsaet,JurisicVidali,Urlep,Vidali} 
and \cite[Section~6.3]{SurveyDRG},  
and thus it may provide some extra information on a possible combinatorial 
structure of the graph (perhaps, it was first observed by Cameron,  
Goethals and Seidel in \cite{SRGSRG}). 
Unfortunately, analysing these equations is rather complicated, especially, 
for a family of distance-regular graphs with unbounded diameter, as the numbers 
of equations and unknowns depend on the diameter.

Many of the Krein parameters vanish when $\Gamma$ is $Q$-polynomial, 
as in this case they satisfy the triangle inequality. This suggests that the triple intersection 
numbers may play an important role in the problem of classification of $Q$-polynomial 
distance-regular graphs. 
In particular, Ivanov and Shpectorov \cite{IVHer} proved that 
a distance-regular graph $\Gamma$ with the same intersection numbers as 
the graph $Her(n,q)$ of Hermitian $n\times n$-forms over $\mathbb{F}_q$ 
is indeed isomorphic to $Her(n,q)$ if $[2,1,1]=0$ holds for any three pairwise 
adjacent vertices $x,y,z$ of $\Gamma$.
Terwilliger \cite[Corollary~2.13]{KiteFree} completed the characterization of $Her(n,q)$ 
by its intersection numbers by observing that $[i,i-1,i-1]=0$ with $2\leq i \leq D$ holds 
for any three pairwise adjacent vertices $x,y,z$ of a distance-regular graph with classical parameters 
$(D,b,\alpha,\beta)$ where $b<-1$.
To do so, Terwilliger \cite[Theorem~2.11]{KiteFree} (cf. Dickie \cite[Theorem~2.1]{Dickie}) 
proved that, for a distance-regular graph with classical parameters $(D,b,\alpha,\beta)$, 
one has $[i,i-1,i-1]=\tau_i[2,1,1]$, where $\tau_i$ is a real scalar 
that depends only on the parameters but not on the particular choice of three 
pairwise adjacent vertices $x,y,z$, and moreover, 
$b<-1$ (which is the case for $Her(n,q)$) implies $\tau_i<0$.

For our purposes, we shall analyse the triple intersection numbers of the type 
$[i,i+1,i+1]$. Theorem \ref{theo-kites} below can be found in \cite[Theorem~3.3]{GavrKoolen} 
in its general form, i.e., not restricted to the case of classical parameters. 

\begin{theo}\label{theo-kites}
Let $\Gamma$ be a distance-regular graph with classical parameters $(D,b,\alpha,\beta)$ 
and diameter $D\geq 3$.  
Suppose that $x,y,z\in \Gamma$ satisfy $x\sim y, x\sim z$ and $\partial(y,z)=j$, $j\in \{1,2\}$.
Then
\begin{align}
[i,i+1,i+1] &= p_{i,i+1}^1\big(\frac{\sigma_i}{b_1}[1,2,2]+\rho_{ij}\big),~~(1\leq i\leq D-1)\label{eq-spear}
\end{align}
where
$\sigma_i = \genfrac{[}{]}{0pt}{}{i}{1}$ (see Eq. (\ref{eq-classical-qbinom})), and
%
\begin{align*}
\rho_{i1} &= - b\genfrac{[}{]}{0pt}{}{i-1}{1},&
\rho_{i2} = - b\genfrac{[}{]}{0pt}{}{i-1}{1} + \frac{b}{b_1}\big(c_i - \sigma_i\big).
\end{align*}
\end{theo}

Note that, for a pair of vertices $y,z\in \Gamma(x)$, we have 
\begin{align}\label{eq-211to111}
[1,2,2] &= 
\left\{
\begin{matrix}
b_1&\mbox{if $y=z,$} &  & \\
b_1-a_1+1+[1,1,1]&\mbox{if $y\sim z,$} &  & \\
b_1-a_1-1+[1,1,1]&\mbox{if $\partial(y,z)=2.$}
\end{matrix}
\right.
\end{align}

\subsection{Local eigenvalues and the Terwilliger polynomial}\label{sect-Tpoly}

We first recall a basic result about the so-called {\it local eigenvalues} 
of a distance-regular graph, i.e., the eigenvalues of its local graphs: 
Theorem \ref{theo-localev} below follows from \cite[Theorems~4.4.3,~4.4.4]{BCN}.

\begin{theo}\label{theo-localev}
Let $\Gamma$ be a distance-regular graph of diameter $D\geq 3$ and with eigenvalues 
$b_0=\theta_0>\theta_1>\ldots>\theta_D$, whose multiplicities are 
$m_0=1,m_1,\ldots,m_D$, respectively. Then, for every vertex $x\in \Gamma$, 
the smallest eigenvalue of the local graph $\Gamma(x)$ is at least 
$\hat{\theta}_1:=-1-\frac{b_1}{\theta_1+1}$. 
If $m_1<b_0$ holds, then the local graph $\Gamma(x)$ has eigenvalue $\hat{\theta}_1$ 
with multiplicity at least $b_0-m_1$.
\end{theo}

Using an algebraic framework behind Theorem \ref{theo-kites}, which is known as 
the Terwilliger (or subconstituent) algebra of a $Q$-polynomial distance-regular graph, 
one can obtain stronger conditions on the eigenvalues of its local graphs. 
Let $\Gamma$ be a distance-regular graph of diameter $D$.
Fix a vertex $x\in \Gamma$, and, for each integer $i$ with $0\leq i\leq D$, 
let $E_i^*:=E_i^*(x)$ 
denote a diagonal matrix with rows and columns indexed 
by the vertex set of $\Gamma$, and defined by 
\[
(E_i^*)_{y,y}=(A_i)_{x,y}
~~~ (y\in \Gamma).
\]

The {\it dual Bose-Mesner} algebra $\mathcal{ M}^*:=\mathcal{ M}^*(x)$ 
with respect to the ({\it base}) vertex $x$ is the matrix algebra 
generated by $E_0^*,E_1^*,\ldots,E_D^*$.
The {\it Terwilliger} (or {\it subconstituent}) algebra 
$\mathcal{ T}:=\mathcal{ T}(x)$ with respect to $x$ is the matrix algebra generated 
by the Bose-Mesner algebra $\mathcal{ M}$ and $\mathcal{ M}^*(x)$, see \cite{SubAlgPaper}. 

Now the triple intersection numbers $[\ell,m,n]_{x,y,z}$ can be expressed in terms of 
the generators of the Terwilliger algebra $\mathcal{T}(x)$ of $\Gamma$ as follows:
\begin{equation}\label{eq-triple-terw}
[\ell,m,n]_{x,y,z}=(E_1^*A_mE_{\ell}^*A_n E_1^*)_{y,z}.
\end{equation} 

We recall that $A_1$ is the adjacency matrix of $\Gamma$, and, 
with an appropriate ordering of the vertices of $\Gamma$, one can see that 
$$
\widetilde{A}:=E_1^*A_1E_1^*=\left(
 \begin{tabular}{ll}
 $N$ & 0\\
 0 & 0
 \end{tabular}
 \right), 
$$
where the principal submatrix $N$ is the adjacency matrix of the local graph $\Gamma(x)$ 
of the base vertex $x$.
With this notation, the equations relating $[i,i+1,i+1]$ to $[1,2,2]$ 
(as in Theorem \ref{theo-kites}, see also \cite[Theorem~2.11]{KiteFree}, \cite[Theorem~2.1]{Dickie})) 
and $[i,i-1,i-1]$ to $[2,1,1]$ (\cite[Theorem~3.3]{GavrKoolen}) yield that
\[
E_1^*A_{i-1}E_i^*A_{i-1}E_1^*\mbox{~~and~~}
E_1^*A_{i}E_{i-1}^*A_{i}E_1^*
\]
are the polynomials of degree at most $2$ in $\widetilde{A}$, $E_1^*$ and $\widetilde{J}:=E_1^*JE_1^*$.
This observation enabled Terwilliger to prove the following strong result 
about the eigenvalues of $\widetilde{A}$, i.e., the eigenvalues of the local graph of $x$. 

\begin{theo}\label{theo-roots}
Let $\Gamma$ be a $Q$-polynomial distance-regular graph with 
classical parameters $(D,b,\alpha,\beta)$, diameter $D\geq 3$ 
{and $|b|\ne 1$}. 
For $i=2,3,\ldots,D-1$, let $T_i(\zeta)$ be a polynomial of degree $4$ defined by 
\[T_i(\zeta):=
-(b^i-1)(b^{i-1}-1)
\times 
\big(\zeta-\beta+\alpha+1\big)
\big(\zeta+1\big)
\big(\zeta+b+1\big)
\big(\zeta-\alpha b\frac{b^{D-1}-1}{b-1}+1\big).
\]
Then, for any vertex of $\Gamma$ and a non-principal eigenvalue $\eta$ of its local graph, 
$T_i(\eta)\geq 0$ holds.
\end{theo}

We call the polynomial $T_i(\zeta)$ the {\it Terwilliger polynomial} of $\Gamma$.
Theorem \ref{theo-roots} was first shown by Terwilliger 
in his ``Lecture note on Terwilliger algebra'' (edited by Suzuki) \cite{SN}. 
The explicit formula of the Terwilliger polynomial was given in our recent paper, 
see \cite[Theorem~4.2,~Proposition~4.3]{GavrKoolen}. 
We refer the reader to \cite{GavrKoolen} for further details, in particular, 
for the general form of Theorem \ref{theo-roots}, which is not restricted to the case 
of classical parameters.

\section{Local graphs of $\Gamma$}\label{sect-prelim}

In this section, we obtain some preliminary results about the local graphs of vertices 
of a distance-regular graph with the same intersection numbers as the Grassmann graph $J_q(2D,D)$. 
Let us first recall some facts about the local structure of the Grassmann graphs, 
see \cite[Chapter~9.3]{BCN} for details.

\begin{reslt}\label{reslt-Jqlocalstr}
\begin{enumerate} 
\item[$(1)$] For every vertex of the Grassmann graph $J_q(n,D)$, 
its local graph is isomorphic to the $q$-clique extension of the 
$\Big(\genfrac{[}{]}{0pt}{}{n-D}{1}_q\times \genfrac{[}{]}{0pt}{}{D}{1}_q\Big)$-grid.
\item[$(2)$] For every pair of vertices at distance $2$ in $J_q(n,D)$, 
their $\mu$-graph is isomorphic to the $(q+1)\times (q+1)$-grid.
\end{enumerate}
\end{reslt}

It follows from Eq. (\ref{eq-gridspectr}) and Results \ref{reslt-Jqlocalstr} and \ref{reslt-cliqueextens} 
that the $q$-clique extension of 
the $(\genfrac{[}{]}{0pt}{}{D}{1}_q\times \genfrac{[}{]}{0pt}{}{D}{1}_q)$-grid, 
which is a local graph in the Grassmann graph $J_q(2D,D)$, 
has spectrum:
\begin{equation}\label{eq-locspectrum2d}
[\hat{\theta}_1]^{g(\hat\theta_1)}, [-1]^{g(-1)}, [\hat\theta_D]^{g(\hat\theta_D)}, [a_1]^1, 
\end{equation}
where the valency $a_1:=q\Big(2\genfrac{[}{]}{0pt}{}{D}{1}_q-1\Big)-1$, and 
\begin{align}
\hat{\theta}_1 &:=-q-1, & \hat{\theta}_D&:=q\Big(\genfrac{[}{]}{0pt}{}{D}{1}_q-1\Big)-1,\\
g(\hat{\theta}_1)&:=\Big(\genfrac{[}{]}{0pt}{}{D}{1}_q-1\Big)^2, & 
g(\hat{\theta}_D)&:=2\Big(\genfrac{[}{]}{0pt}{}{D}{1}_q-1\Big),
\end{align}
\begin{equation}\label{eq-locspectrum2d-2}
g(-1) :=(q-1)\genfrac{[}{]}{0pt}{}{D}{1}_q^2.
\end{equation}

We now formulate the main result of this section.

\begin{prop}\label{prop-localgraphsGamma} 
Let $\Gamma$ be a distance-regular graph 
with classical parameters $(D,q,q,\genfrac{[}{]}{0pt}{}{D+1}{1}_q-1)$ for 
some integers $D\geq 4$ and $q\geq 2$.
The following holds for the local graph $\Delta=\Gamma(x)$ of any vertex $x\in \Gamma$.
\begin{itemize}
\item[$(1)$] $|\Delta(y,z)|\equiv q-2~(\mathsf{mod~}\genfrac{[}{]}{0pt}{}{D-1}{1}_q)$ 
for any pair $y,z$ of vertices of $\Delta$ with $y\sim z$.
\item[$(2)$] $|\Delta(y,z)|=2q$ for any pair $y,z$ of distinct vertices of $\Delta$ 
with $y\not\sim z$. 
\item[$(3)$] $\Delta$ is cospectral to the $q$-clique extension of 
the $\Big(\genfrac{[}{]}{0pt}{}{D}{1}_q\times \genfrac{[}{]}{0pt}{}{D}{1}_q\Big)$-grid.
\end{itemize}
\end{prop}

Clearly, the local graphs in $J_q(2D,D)$ 
satisfy the conclusion of the proposition, 
and our job in Section \ref{sect-main} will be to prove 
that a graph $\Delta$ satisfying Statements (1)--(3) of Proposition \ref{prop-localgraphsGamma}
is indeed isomorphic to the local graphs in $J_q(2D,D)$, 
i.e., the $q$-clique extension of the $\Big(\genfrac{[}{]}{0pt}{}{D}{1}_q\times \genfrac{[}{]}{0pt}{}{D}{1}_q\Big)$-grid.

The proof of Proposition \ref{prop-localgraphsGamma} is given by Lemma \ref{lemma-147} 
and Proposition \ref{prop-locspectra} below, and it exploits the $Q$-polynomial property of $\Gamma$, 
which makes possible to analyse its triple intersection numbers (see Section \ref{sect-triple})
and to apply the Terwilliger polynomial (see Section \ref{sect-Tpoly}).

\begin{lemma}\label{lemma-ddd}
Let $\Gamma$ be a distance-regular graph with 
classical parameters $(D,q,q,\genfrac{[}{]}{0pt}{}{n-D+1}{1}_q-1)$ 
for some integers $n\geq 2D$, $D\geq 3$ and $q\geq 2$. 
Suppose that $x,y,z\in \Gamma$ satisfy $x\sim y, x\sim z$. 
Then 
\begin{align*}
[D-1,D,D]_{x,y,z} &= 
\left\{
\begin{matrix}
\gamma
\Big([1,1,1]_{x,y,z}+q^3\genfrac{[}{]}{0pt}{}{n-D-1}{1}
\big(\genfrac{[}{]}{0pt}{}{D-1}{1}+1\big)
-q\genfrac{[}{]}{0pt}{}{n-D}{1}\genfrac{[}{]}{0pt}{}{D}{1}+2\Big) & \mbox{if $y\sim z,$} \\\\
\gamma
\Big([1,1,1]_{x,y,z}+q^2(q^{n-D-1}-1)-q(q^{D-1}+1)\Big) &\mbox{if $\partial(y,z)=2,$}
\end{matrix}
\right.
\end{align*}
where $\gamma:=q^{D^2-4}\frac{\genfrac{[}{]}{0pt}{}{n-D-1}{D-1}}{\genfrac{[}{]}{0pt}{}{n-D-1}{1}}$.
\end{lemma}
\proof Substituting Eq. (\ref{eq-211to111}) and the classical parameters from 
the statement of the lemma 
into Eq. (\ref{eq-spear}) shows the result.
\wbull

The following lemma implies Statements (1) and (2) of Proposition \ref{prop-localgraphsGamma}.

\begin{lemma}\label{lemma-147}
Let $n=2D$ hold. Then, with the notation as in Lemma $\ref{lemma-ddd}$, one has 
 \begin{align*}
[1,1,1]_{x,y,z} &\equiv~ 
\left\{
\begin{matrix}
q-2~(\mathsf{mod~}\genfrac{[}{]}{0pt}{}{D-1}{1}) &\mbox{if $y\sim z,$} &  & \\\\
2q~(\mathsf{mod~}\genfrac{[}{]}{0pt}{}{D-1}{1})  &\mbox{if $\partial(y,z)=2.$}
\end{matrix}
\right.
\end{align*}
In particular, if $D\geq 4$, then every $\mu$-graph in $\Gamma$ is regular with valency $2q$.
\end{lemma}
\proof Suppose that $n=2D$. 
If $\partial(y,z)=1$, then it follows from Lemma \ref{lemma-ddd} that 
\begin{align*}
[D-1,D,D]_{x,y,z} &= {\displaystyle q^{D^2-4}\frac{1}{\genfrac{[}{]}{0pt}{}{D-1}{1}}
\Big([1,1,1]_{x,y,z}+q^3\genfrac{[}{]}{0pt}{}{D-1}{1}
\big(\genfrac{[}{]}{0pt}{}{D-1}{1}+1\big)
-q\genfrac{[}{]}{0pt}{}{D}{1}^2+2\Big),
}
\end{align*}
where we observe that $\mathsf{gcd}(q^{D^2-4},\genfrac{[}{]}{0pt}{}{D-1}{1})=1$.
Therefore, $\genfrac{[}{]}{0pt}{}{D-1}{1}$ divides 
\[
[1,1,1]_{x,y,z}-q\genfrac{[}{]}{0pt}{}{D}{1}^2+2 = 
[1,1,1]_{x,y,z}-q\Big(\genfrac{[}{]}{0pt}{}{D}{1}^2-1\Big)+2-q =
[1,1,1]_{x,y,z}-q^2\genfrac{[}{]}{0pt}{}{D-1}{1}\Big(\genfrac{[}{]}{0pt}{}{D}{1}+1\Big)+2-q
\]
and thus $[1,1,1]_{x,y,z}\equiv q-2~(\mathsf{mod~}\genfrac{[}{]}{0pt}{}{D-1}{1})$.

Similarly, if $\partial(y,z)=2$, then it follows from Lemma \ref{lemma-ddd} that 
\begin{align*}
[D-1,D,D]_{x,y,z} &= {\displaystyle q^{D^2-4}
\frac{1}{\genfrac{[}{]}{0pt}{}{D-1}{1}}
\Big( 
[1,1,1]_{x,y,z}+q^2(q^{D-1}-1)-q(q^{D-1}+1)
\Big),
}
\end{align*}
where again $\mathsf{gcd}(q^{D^2-4},\genfrac{[}{]}{0pt}{}{D-1}{1})=1$, 
and $\genfrac{[}{]}{0pt}{}{D-1}{1}$ divides $q^2(q^{D-1}-1)$.
Therefore, $\genfrac{[}{]}{0pt}{}{D-1}{1}$ divides 
\[
[1,1,1]_{x,y,z}-q(q^{D-1}+1)=[1,1,1]_{x,y,z}-2q-q(q^{D-1}-1),
\]
and thus $[1,1,1]_{x,y,z}\equiv 2q~(\mathsf{mod~}\genfrac{[}{]}{0pt}{}{D-1}{1})$.

The $\mu$-graph of $y,z$ in $\Gamma$ contains precisely $c_2=(q+1)^2$ vertices,
and $[1,1,1]_{x,y,z}$ is the valency of $x$ in the $\mu$-graph $\Gamma(y,z)$. 
If $D\geq 4$, then $\genfrac{[}{]}{0pt}{}{D-1}{1}>q^{D-2}\geq q^2$, 
and hence $[1,1,1]_{x,y,z}=2q$. 
\wbull

The following proposition proves Statement (3) of Proposition \ref{prop-localgraphsGamma}.

\begin{prop}\label{prop-locspectra}
Let $\Gamma$ be a distance-regular graph with 
classical parameters $(D,q,q,\genfrac{[}{]}{0pt}{}{D+1}{1}_q-1)$ for 
some integers $D\geq 3$ and $q\geq 2$.
Then, for every vertex $x\in\Gamma$, its local graph $\Gamma(x)$ is 
cospectral to the $q$-clique extension of 
the $\Big(\genfrac{[}{]}{0pt}{}{D}{1}_q\times \genfrac{[}{]}{0pt}{}{D}{1}_q\Big)$-grid. 
\end{prop}
\proof 
It follows from 
Theorem \ref{theo-roots} that all Terwilliger polynomials $T_i(\zeta)$, $2\leq i\leq D-1$, 
of $\Gamma$ have the following roots:
\[
-q-1 ~< ~ -1 ~<~
q^2\genfrac{[}{]}{0pt}{}{D-1}{1}_q-1 ~\leq~ \genfrac{[}{]}{0pt}{}{D+1}{1}_q-q-2,
\]
while their leading term coefficients are negative, and, moreover, the two largest roots coincide:
\[
q^2\genfrac{[}{]}{0pt}{}{D-1}{1}_q-1=\genfrac{[}{]}{0pt}{}{D+1}{1}_q-q-2.
\]

Hence, by Theorem \ref{theo-roots}, a non-principal eigenvalue $\eta$ of the local graph $\Gamma(x)$
satisfies:
\begin{equation}\label{Eq-ev-bounds}
-q-1\leq \eta\leq -1\mbox{~~or~~}
\eta=\hat{\theta}_D:=q^2\genfrac{[}{]}{0pt}{}{D-1}{1}_q-1.
\end{equation}

Further, Result \ref{reslt-classparam} implies that
\begin{align*}
\hat{\theta}_1 &:= -1-\frac{b_1}{\theta_1+1}=-1-
\frac{q^3\genfrac{[}{]}{0pt}{}{D-1}{1}_q\genfrac{[}{]}{0pt}{}{D-1}{1}_q}
{q^{2}\genfrac{[}{]}{0pt}{}{D-1}{1}_q\genfrac{[}{]}{0pt}{}{D-1}{1}_q-
\genfrac{[}{]}{0pt}{}{1}{1}_q+1}=-q-1,\\
b_0-m_1 &= q\genfrac{[}{]}{0pt}{}{D}{1}_q\genfrac{[}{]}{0pt}{}{D}{1}_q
-
\genfrac{[}{]}{0pt}{}{2D}{1}_q + \genfrac{[}{]}{0pt}{}{2D}{0}_q=
\Big(\genfrac{[}{]}{0pt}{}{D}{1}_q-1\Big)^2>0,
\end{align*}
and hence, by Theorem \ref{theo-localev}, the local graph $\Gamma(x)$ has eigenvalue $\hat{\theta}_1=-q-1$ 
with multiplicity $h(\hat{\theta}_1)$ at least $b_0-m_1$.

We observe that the valency $a_1$ of the local graph $\Gamma(x)$ cannot be a non-principal eigenvalue 
of $\Gamma(x)$, 
as it does not satisfy Eq. (\ref{Eq-ev-bounds}), and therefore $\Gamma(x)$ is connected and has spectrum:
\[
[\hat{\theta}_1]^{h(\hat\theta_1)}, \eta_1, \ldots, \eta_s,
[-1]^{h(-1)}, [\hat{\theta}_D]^{h(\hat{\theta}_D)},
[a_1]^1, 
\]
where $h(\hat{\theta}_1)\geq b_0-m_1=(\genfrac{[}{]}{0pt}{}{D}{1}_q-1)^2$, 
and 
$s=b_0-h(\hat{\theta}_1)-h(-1)-h(\hat{\theta}_D)-1$ is the number of 
eigenvalues $\eta_i$, $1\leq i\leq s$, satisfying $\hat\theta_1<\eta_i<-1$. 

The spectrum of the $q$-clique extension of 
the $\Big(\genfrac{[}{]}{0pt}{}{D}{1}_q\times \genfrac{[}{]}{0pt}{}{D}{1}_q\Big)$-grid
is given by Eqs. 
(\ref{eq-locspectrum2d})--(\ref{eq-locspectrum2d-2}).
We shall show that $s=0$, i.e., $\{\eta_1,\ldots,\eta_s\}=\emptyset$ and
$h(\eta)=g(\eta)$ for $\eta\in \{\hat{\theta}_1,-1,\hat{\theta}_D\}$.
Let us denote  
\begin{align*}
e_1 &:=h(\hat{\theta}_1)-g(\hat{\theta}_1)~~
(\mbox{so~that~}e_1\geq 0),\\
e_{-1} &:=h(-1)-g(-1),\\
e_D &:=h(\hat{\theta}_D)-g(\hat{\theta}_D).
\end{align*}

Applying Eqs. (\ref{eq-trA-1}) and (\ref{eq-trA-2}) to the adjacency matrices 
of $\Gamma(x)$ and the $q$-clique extension of 
the $\Big(\genfrac{[}{]}{0pt}{}{D}{1}_q\times \genfrac{[}{]}{0pt}{}{D}{1}_q\Big)$-grid, we obtain:
\begin{align*}
1+h(\hat\theta_D)+h(-1)+s+h(\hat\theta_D) &= 
1+g(\hat\theta_D)+g(-1)+g(\hat\theta_1) = b_0,
\\
a_1+\hat\theta_D h(\hat\theta_D)-h(-1)+\sum_{i=0}^s \eta_i+\hat\theta_1 h(\hat\theta_1) &= 
a_1+\hat\theta_D g(\hat\theta_D)-g(-1)+\hat\theta_1 g(\hat\theta_1) = 0,
\\
a_1^2+\hat\theta_D^2h(\hat\theta_D)+h(-1)+\sum_{i=0}^s \eta_i^2+\hat\theta_1^2h(\hat\theta_1) &= 
a_1^2+\hat\theta_D^2g(\hat\theta_D)+g(-1)+\hat\theta_1^2g(\hat\theta_1) = b_0a_1,
\end{align*}
which gives
\begin{align}
e_D+e_{-1}+s+e_{1} & =0,\label{eq-difftr0}
\\
\hat\theta_D e_{D}-e_{-1}+\sum_{i=1}^s \eta_i+\hat\theta_1 e_{1} &= 0,\label{eq-difftr1}
\\
\hat\theta_D^2 e_{D}+e_{-1}+\sum_{i=1}^s \eta_i^2+\hat\theta_1^2 e_{1} &= 0.\label{eq-difftr2}
\end{align}

Multiplying Eq. (\ref{eq-difftr0}) by $\hat{\theta}_D$, Eq. (\ref{eq-difftr1}) 
by $\hat{\theta}_D-1$, and subtracting their sum of Eq. (\ref{eq-difftr2}) gives:
\[
\sum_{i=1}^s (\eta_i+1)(\eta_i-\hat\theta_D)+e_{1}(\hat\theta_1+1)(\hat\theta_1-\hat\theta_D)=0,
\]
which forces $s=e_1=0$, as $e_1\geq 0$, $(\hat\theta_1+1)(\hat\theta_1-\hat\theta_D)>0$, 
and $(\eta_i+1)(\eta_i-\hat\theta_D)>0$ for any $\eta_i$ satisfying $\hat{\theta}_1<\eta_i<-1$.

Thus, Eqs. (\ref{eq-difftr0}) and (\ref{eq-difftr1}) become:
\begin{align*}
e_D+e_{-1} & =0,
\\
\hat\theta_D e_{D}-e_{-1} &= 0,
\end{align*}
which shows $e_{D}=e_{-1}=0$. This proves the proposition. \wbull

\section{Main result}\label{sect-main}

In Section \ref{sect-bigcliques}, we prove that the local graphs of a distance-regular 
graph $\Gamma$ satisfying the conditions of Proposition \ref{prop-localgraphsGamma} are 
indeed isomorphic to the $q$-clique extension of a square grid if $D$ is large enough. 
In Section \ref{sect-local}, we recall a theorem by Numata, Cohen and Cooperstein, 
and show that applying it to $\Gamma$ completes the proof of Theorem \ref{theo-main-1}.

\subsection{Spectral characterization of the local graphs}\label{sect-bigcliques}


\begin{prop}\label{prop-main} 
Let $\Delta$ be a graph satisfying the following conditions 
for some $q,D\in \mathbb{N}$, $q\geq 2$.
\begin{itemize}
\item[$(1)$] $|\Delta(y,z)|\equiv q-2~(\mathsf{mod~}\genfrac{[}{]}{0pt}{}{D-1}{1}_q)$ 
for any pair $y,z$ of vertices of $\Delta$ with $y\sim z$.
\item[$(2)$] $|\Delta(y,z)|=2q$ for any pair $y,z$ of distinct vertices of $\Delta$ 
with $y\not\sim z$. 
\item[$(3)$] $\Delta$ is cospectral to the $q$-clique extension of 
the $\Big(\genfrac{[}{]}{0pt}{}{D}{1}_q\times \genfrac{[}{]}{0pt}{}{D}{1}_q\Big)$-grid.
\end{itemize}

If $D\geq \chi(q)$ (see Eq. (\ref{eq-chi})), 
then $\Delta$ is isomorphic to the $q$-clique extension of 
the $\Big(\genfrac{[}{]}{0pt}{}{D}{1}_q\times \genfrac{[}{]}{0pt}{}{D}{1}_q\Big)$-grid.
\end{prop}

The proof of Proposition \ref{prop-main} is based on the idea from \cite{Riaz}.
For the rest of this section, let $\Delta$ be a graph satisfying the condition of Proposition 
\ref{prop-main}. 
To simplify the notation, put $r:=\genfrac{[}{]}{0pt}{}{D}{1}_q$ and 
$k:=q(2r-1)-1$ (note that $k$ is the valency of a vertex of $\Delta$).
We call a maximal clique of $\Delta$ a {\it line} if it contains at least
$\kappa qr+1$ vertices where $\kappa$ is any real number satisfying 
$\frac{2}{3}+\frac{5q-4}{3qr}<\kappa\leq 1-\frac{1}{qr}$. 
We first show in Lemma \ref{lemma-2lines} that every vertex of $\Delta$ lies in exactly 
two lines. In Lemma \ref{lemma-final}, we then prove that every line has the same size $qr$ 
and every two non-trivially intersecting lines share precisely $q$ vertices, 
which reveals the structure of $\Delta$.

Fix a vertex $\infty$ of $\Delta$, and let the vertices of $\Delta(\infty)$ have 
the valencies $\lambda_1,\ldots,\lambda_k$ in $\Delta(\infty)$. 

\begin{lemma} \label{lemma-localpreliminaries}
\begin{itemize}
\item[$(1)$] A clique of $\Delta$ has size at most $qr$. 
\item[$(2)$] For any vertex of $\Delta$, its local graph contains a coclique of size at most $(q+1)^2$.
\item[$(3)$] The following equalities hold:
\begin{align}
\sum_{i=1}^k \lambda_i &= q^2(2r^2 - 1) - 3q(2r-1) + 2,\label{eq-local-valencies-sum} \\
\sum_{i=1}^k \lambda_i^2 &= q^3 (2 r^3 + 2 r^2 - 4 r + 1) + q^2 (-12 r^2 + 4 r + 3) + 8 q (2 r - 1) - 4. \label{eq-local-valencies-square-sum}
\end{align}
\end{itemize}
\end{lemma}
\proof We recall that the spectrum of $\Delta$ is given by 
Eqs. (\ref{eq-locspectrum2d})--(\ref{eq-locspectrum2d-2}):
\[
[-q-1]^{(r-1)^2}, [-1]^{(q-1)r^2}, [q(r-1)-1]^{2(r-1)}, [k]^1.
\]

Let $L$ be a clique of size $\ell$ in $\Delta$.
The partition $\{L,\Delta\setminus L\}$ 
of the vertex set of $\Delta$ has quotient matrix
\begin{equation*}
Q=\left(
  \begin{array}{cc}
    \ell-1 & q(2r-1)-\ell \\
    \frac{(q(2r-1)-\ell)\ell}{qr^{2}-\ell} & q(2r-1)-1-\frac{(q(2r-1)-\ell)\ell}{qr^{2}-\ell}\\
  \end{array}
\right)
\end{equation*}
with eigenvalues $k=q(2r-1)-1$ and $\ell-1-\frac{(q(2r-1)-\ell)\ell}{qr^{2}-\ell}$. 
By Result \ref{lema-matrix-interlacing},
we obtain that the second eigenvalue of the quotient matrix $Q$ is at most $q(r-1)-1$, i.e., 
\begin{equation*}
\ell-1-\frac{(q(2r-1)-\ell)\ell}{qr^{2}-\ell}\leq q(r-1)-1
\end{equation*}
holds, which simplifies to $\ell\leq qr$.
This shows (1).

The smallest eigenvalue of a complete bipartite graph with parts of size $1$ and $c$ 
is $-\sqrt{c}$. By Result \ref{lem-interlacing}, 
a $c$-coclique may exist in a local graph of $\Delta$ 
if $-\sqrt{c}\geq -q-1$. This shows (2).

The number of triangles through $\infty$ equals $\frac{1}{2}\sum_{i=1}^k\lambda_i$. On the other hand, as $\Delta$ has just 4 distinct eigenvalues, it follows by Result \ref{reslt-vanDam}(1) that 
\begin{align*}
\frac{1}{2}\sum_{i=1}^k\lambda_i 
&= q^2r^2 - q^2/2 - 3qr + 3q/2 + 1,
\end{align*}
which shows Eq. (\ref{eq-local-valencies-sum}).

By Condition (2) of Proposition \ref{prop-main}, $\Delta$ contains precisely $(qr^2-k-1){2q\choose 2}$ 
quadrangles through $\infty$ that do not have a diagonal edge incident to $\infty$.
The number of quadrangles having a diagonal edge incident to $\infty$ equals 
$\sum_{i=1}^k {\lambda_i \choose 2}$. The total number of quadrangles through $\infty$ 
is given by Result  \ref{reslt-vanDam}(2). Combining these facts gives 
Eq. (\ref{eq-local-valencies-square-sum}) and completes the proof of the lemma.
\wbull

Further, combining Eqs. (\ref{eq-local-valencies-sum}) and (\ref{eq-local-valencies-square-sum}), 
we obtain the following useful equation:
\begin{equation}\label{eq-useful}
\sum\limits_{i=1}^{k}\big(\lambda_i-(qr-2)\big)^{2}=q^2(r-1)^2(q-1).
\end{equation}


%
%
%

\begin{lemma}\label{lemma-2lines}
If $r\geq \frac{q^{\chi(q)}-1}{q-1}$ 
holds, 
then each vertex in $\Delta$ lies in exactly two lines.
\end{lemma}
\proof 
Let $C$ be a maximal coclique of $\Delta(\infty)$ with vertex set $\{x_1,x_2,\ldots,x_c\}$. 
By Lemma \ref{lemma-localpreliminaries}(2), one has $c:=|C|\leq (q+1)^2$.
We define
\[
P:=\{y\in \Delta(\infty)\mid y \mbox{ has at least two neighbours in }C\},
\]
\[
U_i:=\{x_i\}\cup \{y\in \Delta(\infty)\mid y \mbox{ has only }x_i\mbox{ as its neighbour in }C\}~~
(1\leq i\leq c).
\]

The maximality of $C$ implies that $\{P,U_1,\ldots,U_c\}$ is a partition of the vertex set of $\Delta(\infty)$ and each $U_i$ induces a complete subgraph in $\Delta(\infty)$.
Put $p:=|P|$, and $u_i:=|U_i|$ for $1\leq i\leq c$.

By Lemma \ref{lemma-localpreliminaries}(1), it follows that $u_i\leq qr-1$.
As $x_i$ and $x_j$ with $i\ne j$ 
have at most $2q-1$ common neighbours in $\Delta(\infty)$, we obtain
\begin{equation}\label{b-eqn}
p\leq (2q-1){c\choose 2}<qc(c-1)\leq q^2(q+1)^2(q+2).
\end{equation}

Let $t$ denote the number of edges in $\Delta(\infty)$. 
Then $2t$ equals $\sum_{i=1}^k\lambda_i$, which is given 
by Eq. (\ref{eq-local-valencies-sum}). 
On the other hand, we observe that each $U_i$ contains $u_i(u_i-1)/2$ edges, 
$P$ contains at most $p(p-1)/2$ edges, and there are at most $p(k-p)$ edges 
between $P$ and $\cup_{i=1}^c U_i$, while there are at most $(c-1)(2q-1)u_i$ edges 
between $U_i$ and $\cup_{j\ne i}U_j\cup P$. Thus, we obtain:
\[
2t\leq \sum_{i=1}^c u_i(u_i-1)+p(p-1)+
2p(k-p)+\big(\sum_{i=1}^c(c-1)(2q-1)u_i-p(k-p)\big),
\]
which, by using Eq. (\ref{b-eqn}) and $\sum_{i=1}^cu_i=k-p<k$, simplifies to 
\begin{equation}\label{eq-eps-bound}
2t <\sum_{i=1}^c u_i^2+2rq^3(q+2)((q+1)^2+2).
\end{equation}


Assume that there exists at most one line in $\Delta$ through $\infty$.
Then, by Lemma \ref{lemma-localpreliminaries}(1), for at most one $i^*\in \{1,2,\ldots,c\}$ 
we have $u_{i^*}\leq qr-1$, so that $u_i\leq u_{i^*}$ and 
$u_i\leq \kappa qr$ for each $i\in\{1,2,\ldots,c\}\setminus \{i^*\}$. 
Since $u_i\geq u_j$ implies that, for any $\epsilon>0$,  
\[
(u_{i}+\epsilon)+(u_j-\epsilon)^2=u_{i}^2+u_j^2+2\epsilon(u_{i}-u_j)+2\epsilon^2>u_{i}^2+u_j^2,
\]
one can see that 
\begin{equation}\label{eq-sum-u-sq}
\sum_{i=1}^c u_i^2\leq (qr-1)^2+(\kappa qr)^2+(q(2r-1)-1-(qr-1)-\kappa qr)^2.
\end{equation}


Combining Eqs. (\ref{eq-eps-bound}) and (\ref{eq-sum-u-sq}) with Eq. (\ref{eq-local-valencies-sum}), 
we obtain the inequality
\begin{align*}
2q^2r^2 - q^2 - 6qr + 3q + 2  &< (qr-1)^2+(\kappa qr)^2+q^2(r(1-\kappa)-1)^2+2rq^3(q+2)((q+1)^2+2),
\end{align*}
which violates if $r\geq \frac{q^{\chi(q)}-1}{q-1}$. 

Thus, if $r\geq \frac{q^{\chi(q)}-1}{q-1}$ holds, we obtain that the vertex $\infty$ lies in at 
least two lines. 
Let $L_1$ and $L_2$ be two such lines. By Condition (2) of Proposition \ref{prop-main}, 
$|L_1\cap L_2|\leq 2q$ holds, and hence $\Delta(\infty)\setminus (L_1\cup L_2)$ contains at most 
\[
k-(2\kappa qr-(2q-1))=2qr+q-2-2\kappa qr<
2qr+q-2-\frac{4qr+10q-8}{3}<\kappa qr\]
vertices. This implies that $\infty$ lies in at most two lines. 
Since $\infty$ was arbitrarily chosen in $\Delta$, 
this shows the lemma for every vertex of $\Delta$.\wbull


\begin{lemma}\label{lemma-final} 
Suppose that $r\geq \frac{q^{\chi(q)}-1}{q-1}$ holds. Then every line in $\Delta$ contains 
precisely $qr$ vertices, while every two intersecting lines have precisely $q$ vertices in common.
\end{lemma}
\proof According to Lemma \ref{lemma-2lines}, 
let $L_1$ and $L_2$ be the two lines containing $\infty$. 
Put $\Delta_0:=\Delta(\infty)\setminus (L_1\cup L_2)$ and $\delta_0:=|\Delta_0|$.
As $|L_1\cap L_2|\leq 2q$ holds by Condition (2) of Proposition \ref{prop-main}, 
one has
\[
\delta_0\leq k-(2\kappa qr-(2q-1)).
\]

We observe that a vertex $y\in \Delta_0$ is adjacent to 
at most $2q-1$ vertices in each $L_i$, $i=1,2$, and 
to at most $\delta_0-1$ other vertices in $\Delta_0$, i.e., 
its valency $\nu_y$ in the local graph $\Delta(\infty)$ satisfies 
\begin{equation}\label{eq-delta0valency}
\nu_y\leq 2(2q-1)+\delta_0-1\leq 3(2q-1)+k-2\kappa qr-1<
\frac{2}{3}qr+\frac{5q-7}{3}.
\end{equation}

This implies that $|\nu_y-(qr-2)|>\frac{1}{3}qr+\frac{5q-1}{3}$ 
and it follows from Eq. (\ref{eq-useful}) that
\begin{align*}
q^2(r-1)^2(q-1)&=\sum_{i=1}^k(\lambda_i-(qr-2))^2 > 
\sum_{y\in \Delta_0}(\nu_y-(qr-2))^2 \\
&> \frac{\delta_0}{9}(qr+5q-1)^2,
\end{align*}
which forces $\delta_0<8q$ if $r\geq \frac{q^5-1}{q-1}$.

Pick vertices $y_0, y_1, y_2$ such that 
\[y_0\in \Delta(\infty)\setminus (L_1\cup L_2),~~
y_1\in (L_1\setminus L_2)\cup (L_2\setminus L_1),~~y_2\in (L_1\cap L_2)\setminus \{\infty\},
\]
and let us estimate their valencies in the local graph $\Delta(\infty)$.

By $\delta_0<8q$ and Eq. (\ref{eq-delta0valency}), we have 
\begin{equation}\label{eq-y0val}
0\leq |\Delta(\infty,y_0)|<2(2q-1)+8q=12q-2.
\end{equation}

By $\delta_0<8q$, $|L_i|\leq qr$ for $i=1,2$, and $|L_1\cap L_2|\leq 2q$, we obtain that 
\begin{equation}\label{eq-y1val}
q(2r-1)-1-8q-(qr-1)-1<|\Delta(\infty,y_1)|<qr-2+8q+2q-1,
\end{equation}
and 
\begin{equation}\label{eq-y2val}
q(2r-1)-2-8q<|\Delta(\infty,y_2)|\leq q(2r-1)-2.
\end{equation}

Recall that $r=\genfrac{[}{]}{0pt}{}{D}{1}_q$, and, by Condition (1) of Proposition \ref{prop-main}, 
the number $|\Delta(\infty,y_i)|$, $i=1,2,3$, satisfies 
\begin{equation}\label{eq-yimod}
|\Delta(\infty,y_i)|\equiv q-2~(\mathsf{mod~}\genfrac{[}{]}{0pt}{}{D-1}{1}_q),
\end{equation}

Therefore Eqs. (\ref{eq-y0val})--(\ref{eq-yimod}) and $D\geq \chi(q)$ force 
$|\Delta(\infty,y_i)|=\ell_i$ for $i=1,2,3$ where
\[
\ell_0=q-2,~~\ell_1=qr-2,~~\ell_2=2(qr-2)-(q-2).
\]

Put $\delta_1:=|(L_1\setminus L_2)\cup (L_2\setminus L_1)|$ and 
$\delta_2:=|(L_1\cap L_2)\setminus \{\infty\}|$. Then $k=\delta_0+\delta_1+\delta_2$, 
Eqs. (\ref{eq-local-valencies-sum}) and (\ref{eq-local-valencies-square-sum}) give 
the following system of equations:
\begin{align*}
\delta_0+\delta_1+\delta_2&=q(2r-1)-1,
\\
\ell_0\delta_0+\ell_1\delta_1+\ell_2\delta_2&=2q^2r^2 - q^2 - 6qr + 3q + 2,
\\
\ell_0^2\delta_0+\ell_1^2\delta_1+\ell_2^2\delta_2&=q^3 (2 r^3 + 2 r^2 - 4 r + 1) + q^2 (-12 r^2 + 4 r + 3) + 8 q (2 r - 1) - 4,
\end{align*}
which has a unique solution
\[
\delta_0=0,~~\delta_1=2q(r-1),~~\delta_2=q-1.
\]

This implies that $|L_1|=|L_2|=qr$, $|L_1\cap L_2|=q$, and the lemma follows.
\wbull

Let us complete the proof of Proposition \ref{prop-main}.
By Lemmas \ref{lemma-2lines} and \ref{lemma-final}, we find that every vertex $x$ of $\Delta$ lies
in exactly two lines of order $qr$, and the two lines through $x$ have precisely $q$ vertices in common.
Define the following equivalence relation $\mathcal{E}$ on the vertex set of $\Delta$:
\[
x\mathcal{E}x'\mbox{~if and only if~}\{x\}\cup \Delta(x)=\{x'\}\cup \Delta(x').
\]

It follows that every equivalence class of $\mathcal{E}$ is the intersection of two lines, 
and the $q$ vertices in the same equivalence class induce a $q$-clique of $\Delta$. 
Define a graph $\underline{\Delta}$ whose vertices are the equivalence classes of $\mathcal{E}$ 
with two classes $C_1$, $C_2$ being adjacent whenever a vertex from $C_1$ is adjacent to a vertex 
from $C_2$. Then $\underline{\Delta}$ is a regular graph with valency $2(r-1)$, 
and $\Delta$ is the $q$-clique extension of $\underline{\Delta}$. 
The spectrum of $\underline{\Delta}$ follows from Result \ref{reslt-cliqueextens}:
\begin{equation*}
[2(r-1)]^{1},[r-2]^{2(r-1)},[-2]^{(r-1)^{2}},
\end{equation*}
and thus, see Eq. (\ref{eq-gridspectr}), $\underline{\Delta}$ is the $(r\times r)$-grid. 
This shows Proposition \ref{prop-main}. \wbull

\subsection{A local characterization of the Grassmann graphs}\label{sect-local}
In this section we recall the Numata-Cohen-Cooperstein theorem \cite[Theorem~9.3.8]{BCN}
(see Theorem \ref{theo-NumataCohen} below), 
which characterizes a class of distance-regular graphs including 
the Grassmann graphs by their local structure. 

Recall that an $s$-\emph{coclique} of a graph is an induced subgraph 
on $s$ vertices but without edges. 
We call an $s$-coclique simply a \emph{coclique} if we do not refer to its cardinality.

\begin{theo}\label{theo-NumataCohen}
Let $\Gamma$ be a finite connected graph such that
\begin{itemize}
\item[$(i)$] for every pair of vertices $x,y\in \Gamma$ with $\partial(x,y)=2$, 
the $\mu$-graph of $x,y$ is a non-degenerate grid, and
\item[$(ii)$] if $x,y,z\in \Gamma$ induce a $3$-coclique, then $\Gamma(x,y,z)$ is a coclique.
\end{itemize}
Then $\Gamma$ is either a clique, or a Johnson graph $J(n,k)$, or the quotient of 
the Johnson graph $J(2k,k)$ obtained by identifying a $k$-set with the image of its complement 
under the identity or an involution in $Sym(2k)$ with at least $10$ fixed points (i.e., a folded 
Johnson graph), or a Grassmann graph $J_q(n,D)$ over a finite field $\mathbb{F}_q$.
\end{theo} 

\begin{coro}\label{coro-final}
Let $\Gamma$ be a distance-regular graph 
with classical parameters $(D,q,q,\genfrac{[}{]}{0pt}{}{n-D+1}{1}_q-1)$ 
for some integers $n\geq 2D$, $D\geq 2$ and $q\geq 2$. 
Suppose that, for every vertex $x\in \Gamma$, its local graph $\Gamma(x)$ 
is isomorphic to the $q$-clique extension of the 
$\Big(\genfrac{[}{]}{0pt}{}{n-D}{1}_q\times \genfrac{[}{]}{0pt}{}{D}{1}_q\Big)$-grid.
Then $q$ is a prime power and $\Gamma$ is isomorphic to $J_q(n,D)$.
\end{coro}
\proof As the $q$-clique extension of the 
$\Big(\genfrac{[}{]}{0pt}{}{n-D}{1}_q\times \genfrac{[}{]}{0pt}{}{D}{1}_q\Big)$-grid 
does not contain a 3-claw (a complete bipartite subgraph with parts of size 1 and 3), 
we see that $\Gamma$ satisfies Condition $(ii)$ of Theorem \ref{theo-NumataCohen}.
We shall show that $\Gamma$ satisfies Condition $(i)$ of Theorem \ref{theo-NumataCohen}.
Let $x,y$ be a pair of vertices of $\Gamma$ with $\partial(x,y)=2$, 
and let $M$ denote their $\mu$-graph, which contains precisely $c_2=(q+1)^2$ vertices. 
We observe that the local graph of any vertex $u\in M$ is the disjoint union of two $q$-cliques 
(consider the $\mu$-graph of $x$ and $y$ in the local graph of $u$). 
Therefore, the edge set of $M$ can be partitioned into the edge sets of all maximal $(q+1)$-cliques, 
while each vertex of $M$ lies in two of these cliques.
By the criterion of Krausz \cite{Krausz}, $M$ is a {\it line} graph, namely, 
the line graph of a complete bipartite graph with parts of size $q+1$, i.e., 
the $(q+1)\times (q+1)$-grid.\wbull 

Theorem \ref{theo-main-1} follows from Propositions \ref{prop-localgraphsGamma} and 
\ref{prop-main} and Corollary \ref{coro-final}. 
Theorem \ref{theo-main} follows from Theorem \ref{theo-main-1} and Result \ref{reslt-classparam}.

{\bf Acknowledgements}
\smallskip


The research of Alexander Gavrilyuk was supported by BK21plus Center for Math Research 
and Education at Pusan National University. 
His work (e.g., Proposition \ref{prop-localgraphsGamma}) was also partially 
supported by the Russian Science Foundation (Grant 14-11-00061-P). 

The research of Jack Koolen was partially supported by the National Natural Science Foundation 
of China (Grants No. 11471009 and No. 11671376).


\begin{thebibliography}{00}


\bibitem{Bichara} Bichara, A., Tallini, G.:  
\newblock On a characterization of Grassmann space representing the $h$-dimensional
subspaces in a projective space.
\newblock {\em Ann. Discrete Math.} {\bf 18}, 113--132 (1983) 

\bibitem{Biondi} Biondi, P.: 
\newblock On finite Grassmann spaces.
\newblock {\em Ann. Discrete Math.} {\bf 37}, 69--73 (1988) 

\bibitem{Buss} Bussemaker, F.C., Neumaier, A.:
\newblock Exceptional graphs with smallest eigenvalue $-2$ and related problems.
\newblock {\em Math. Comput.} {\bf 59}, 583--608 (1992) 

\bibitem{Cohen} Cohen, A.M.: 
\newblock On a theorem of Cooperstein.
\newblock {\em Eur. J. Combin.} {\bf 4}, 107--126 (1983) 

\bibitem{CohenCoop} Cohen, A.M., Cooperstein, B.N.: 
\newblock A characterization of some geometries of Lie type.
\newblock {\em Geometriae Dedicata} {\bf 15}, 73--105 (1983) 

\bibitem{Cooperstein} Cooperstein, B.N.: 
\newblock Some geometries associated with parabolic representations of groups of Lie type.
\newblock {\em Can. J. Math.} {\bf 28}, 1021--1031 (1976) 

\bibitem{TwistG} Bang, S., Fujisaki, T., Koolen, J.H.:
\newblock The spectra of the local graphs of the twisted Grassmann graphs.
\newblock {\em Eur. J. Combin.} {\bf 30(3)}, 638--654 (2009) 

\bibitem{BI} Bannai, E., Ito, T.:
\newblock Algebraic combinatorics. I. Association schemes.
\newblock The Benjamin/Cummings Publishing Co., Inc., Menlo Park, CA (1984)

\bibitem{BBI} Bannai, E., Bannai, E., Ito, T.:
\newblock Introduction to Algebraic Combinatorics. (in Japanese)
\newblock Kyoritsu-Shuppan (2016)

\bibitem{biggs} Biggs, N.:
\newblock Algebraic Graph Theory.
\newblock Cambridge University Press, Cambridge (1993)

\bibitem{BCN} Brouwer, A.E., Cohen, A.M., Neumaier, A.:
\newblock Distance-regular graphs.
\newblock {\em Ergebnisse der Mathematik und ihrer Grenzgebiete}, (3), 18. 
Springer-Verlag, Berlin (1989)

\bibitem{BH} Brouwer, A.E., Haemers, W.H.:
\newblock Spectra of Graphs.
\newblock Springer, Heidelberg (2012)

\bibitem{BWil} Brouwer, A.E., Wilbrink, H.A.:
\newblock The structure of near polygons with quads.
\newblock {\em Geometriae Dedicata} {\bf 14(2)}, 145--176 (1983)

\bibitem{SRGSRG} Cameron, P.J., Goethals, J.-M., Seidel, J.J.:
\newblock Strongly regular graphs having strongly regular subconstituents.
\newblock {\em J. Algebra} {\bf 55}, 257--280 (1978)

\bibitem{Cartan} Cartan, \'{E}.:
\newblock Sur une classe remarquable d'espaces de Riemann, I.
\newblock {\em Bulletin de la Soci\'{e}t\'{e} Math\'{e}matique de France} {\bf 54}, 214--216 (1926)

\bibitem{JurisicCoolsaet} Coolsaet, K., Juri\v{s}i\'c, A.:  
\newblock Using equality in the Krein conditions to prove nonexistence of certain distance-regular graphs.
\newblock {\em J. Comb. Theory A} {\bf 115(6)}, 1086--1095 (2008)

\bibitem{CuypersJq} Cuypers, H.:
\newblock The dual of Pasch's axiom.
\newblock {\em Eur. J. Combin.} {\bf 13(1)}, 15--31 (1992)

\bibitem{CuypersBF} Cuypers, H.:
\newblock Two remarks on Huang's characterization of the bilinear forms graphs.
\newblock {\em Eur. J. Combin.} {\bf 13(1)}, 33--37 (1992)

\bibitem{vD95} Van Dam, E.:
\newblock Regular graphs with four eigenvalues.
\newblock {\em Linear Algebra Appl.} {\bf 226-228}, 139--162 (1995)

\bibitem{vDK} Van Dam, E.R., Koolen, J.H.:
\newblock A new family of distance-regular graphs with unbounded diameter.
\newblock {\em Invent. Math.} {\bf 162}, 189--193 (2005) 

\bibitem{SurveyDRG} Van Dam, E.R., Koolen, J.H., Tanaka, H.:
\newblock Distance-regular graphs.
\newblock {\em Electron. J. Comb.}, Dynamic Survey DS22. 

\bibitem{Delsarte} Delsarte, P.:
\newblock An algebraic approach to the association schemes of coding theory.
\newblock {\em Philips Res. Reports Suppl.} {\bf 10} (1973)

\bibitem{Dickie} Dickie, G.:
\newblock Twice $Q$-polynomial distance-regular graphs.
\newblock {\em J. Comb. Theory B} {\bf 68(1)}, 161--166 (1996)

\bibitem{Egawa} Egawa, Y.:
\newblock Characterization of $H(n,q)$ by the parameters.
\newblock {\em J. Comb. Theory A} {\bf 31}, 108--125 (1981)


\bibitem{TwistG2} Fujisaki, T., Koolen, J.H., Tagami, M.: 
\newblock Some properties of the twisted Grassmann graphs.
\newblock {\em Innov. Incidence Geom.} {\bf 3}, 81--87 (2006) 

\bibitem{GavrKoolen} Gavrilyuk, A.L., Koolen, J.H.:
\newblock The Terwilliger polynomial of a $Q$-polynomial distance-regular graph and 
its application to pseudo-partition graphs.
\newblock {\em Linear Algebra Appl.} {\bf 466(1)}, 117--140 (2015)

\bibitem{Bil2d} Gavrilyuk, A.L., Koolen, J.H.:
\newblock A characterization of the graphs of $(d\times d)$-bilinear forms over $\mathbb{F}_2$.
\newblock {\em Combinatorica} (to appear), \href{https://arxiv.org/abs/1511.09435}{arXiv:1511.09435}

\bibitem{Hoffman} Hoffman, A.J.:
\newblock On the polynomial of a graph.
\newblock {\em Am. Math. Mon.} {\bf 70}, 30--36 (1963)

\bibitem{Huang} Huang, T.:
\newblock A characterization of the association schemes of bilinear forms.
\newblock {\em Eur. J. Combin.} {\bf 8}, 159--173 (1987)

\bibitem{Jung} Jungnickel, D.:
\newblock Characterizing geometric designs, II.
\newblock {\em J. Comb. Theory A} {\bf 118}, 623--633 (2011)

\bibitem{IVA2d} Ivanov, A.A., Shpectorov, S.V.:
\newblock The association schemes of dual polar spaces of type $A_{2d-1}(p^f)$ are 
characterized by their parameters if $d\geq 3$.
\newblock {\em Linear Algebra Appl.} {\bf 114-115}, 133--139 (1989)

\bibitem{IVHer} Ivanov, A.A., Shpectorov, S.V.:
\newblock A characterization of the association schemes of Hermitian forms.
\newblock {\em J. Math. Soc. Jpn.} {\bf 43}, 25--48 (1991)

\bibitem{JurisicVidali} Juri\v{s}i\'c, A., Vidali, J.: 
\newblock Extremal 1-codes in distance-regular graphs of diameter 3.
\newblock {\em Design. Code. Cryptogr.} {\bf 65(1-2)}, 29--47 (2012)

\bibitem{Krausz} Krausz, J.:
\newblock D\'{e}monstration nouvelle d'une th\'{e}or\`{e}me de Whitney sur les r\'{e}seaux.
\newblock {\em Mat. Fiz. Lapok} {\bf 50}, 75--85 (1943)

\bibitem{Lore} Lo Re, P.M., Olanda, D.:
\newblock Grassmann spaces.
\newblock {\em J. Geom.} {\bf 17}, 50--60 (1981) 

\bibitem{Martin} Martin, W., Tanaka, H.:
\newblock Commutative association schemes.
\newblock {\em Eur. J. Combin.} {\bf 30(6)}, 1497--1525 (2009) 

\bibitem{Melone} Melone, N., Olanda, D.:
\newblock A characteristic property of the Grassmann manifold representing 
the lines of a projective space.
\newblock {\em Eur. J. Combin.} {\bf 5}, 323--330 (1984) 

\bibitem{Metsch95} Metsch, K.:
\newblock A characterization of Grassmann graphs.
\newblock {\em Eur. J. Combin.} {\bf 16}, 639--644 (1995) 

\bibitem{Metsch971} Metsch, K.:
\newblock Characterization of the folded Johnson graphs of small diameter 
by their intersection arrays.
\newblock {\em Eur. J. Combin.} {\bf 18}, 901--913 (1997) 

\bibitem{Metsch91} Metsch, K.:
\newblock Improvement of Bruck's completion theorem.
\newblock {\em Design. Code. Cryptogr.} {\bf 1}, 99--116 (1991)

\bibitem{Metsch99} Metsch, K.:
\newblock On a characterization of bilinear forms graphs.
\newblock {\em Eur. J. Combin.} {\bf 20}, 293--306 (1999)

\bibitem{Metsch03} Metsch, K.:
\newblock On the characterization of the folded halved cubes by their intersection arrays.
\newblock {\em Design. Code. Cryptogr.} {\bf 29}, 215--225 (2003) 

\bibitem{Metsch97} Metsch, K.:
\newblock On the characterization of the folded Johnson graphs and the folded
halved cubes by their intersection arrays.
\newblock {\em Eur. J. Combin.} {\bf 18}, 65--74 (1997) 

\bibitem{TwistGr3} Munemasa, A.:
\newblock Godsil-McKay switching and twisted Grassmann graphs.
\newblock {\em Design. Code. Cryptogr.} {\bf 84(1-2)}, 173-179 (2017)

\bibitem{TwistGr4} Munemasa, A., Tonchev, V.D.:
\newblock The twisted Grassmann graph is the block graph of a design.
\newblock {\em Innov. Incidence Geom.} {\bf 12}, 1--6 (2011)

\bibitem{Neumaier} Neumaier, A.:
\newblock Characterization of a class of distance regular graphs.
\newblock {\em J. Reine Angew. Math.} {\bf 357}, 182--192 (1985)

\bibitem{Numata} Numata, M.:
\newblock On the graphical characterization of the projective space over a finite field.
\newblock {\em J. Comb. Theory B} {\bf 38}, 143--155 (1985)

\bibitem{RCS} Ray-Chaudhuri, D.K., Sprague, A.P.:
\newblock Characterization of projective incidence structures.
\newblock {\em Geometriae Dedicata} {\bf 5}, 361--376 (1976)

\bibitem{Riaz} Hayat, S., Koolen, J.H., Riaz, M.:
\newblock A spectral characterization of the $s$-clique 
extension of the square grid graphs,
\newblock {\em submitted}.

\bibitem{Shrikhande} Shrikhande, S.S.:
\newblock The uniqueness of the $L_2$ association scheme.
\newblock {\em Ann. Math. Statist.} {\bf 30}, 781--798 (1959)

\bibitem{Shult} Shult, E.E.:
\newblock Characterizations of the Lie incidence geometries.
\newblock In: Lloyd, E.K. (ed.) {\em Surveys in Combinatorics},  
London Mathematical Society Lecture Notes Series {\bf 82}, pp. 157--184. Cambridge University Press (1982)

\bibitem{Shult2} Shult, E.E.:
\newblock A remark on Grassmann spaces and half-spin geometries.
\newblock {\em Eur. J. Combin.} {\bf 15}, 47--52 (1994)

\bibitem{Sprague} Sprague, A.P.:
\newblock Characterization of projective graphs.
\newblock {\em J. Comb. Theory B} {\bf 24}, 294--300 (1978)

\bibitem{Tallini} Tallini, G.:
\newblock On a characterization of the Grassmann manifold representing the lines in a projective space.
\newblock In: Cameron, P.J., Hirschfeld, J.W.P., Hughes, D. (eds.)
{\em Finite Geometries and Designs}, 
London Mathematical Society Lecture Notes Series {\bf 49}, pp. 354--358. Cambridge University Press (1981)

\bibitem{TerwilligerJohnson} Terwilliger, P.:
\newblock The Johnson graph $J(d,r)$ is unique if $(d, r)\ne (2, 8)$.
\newblock {\em Discrete Math.} {\bf 58}, 175--189 (1986) 

\bibitem{SN} Terwilliger, P.:
\newblock Lecture note on Terwilliger algebra (edited by H. Suzuki) (1993)

\bibitem{SubAlgPaper} Terwilliger, P.:
\newblock The subconstituent algebra of an association scheme, I.
\newblock {\em J. Algebr. Comb.} {\bf 1(4)}, 363--388 (1992) 

\bibitem{KiteFree} Terwilliger, P.:
\newblock Kite-free distance-regular graphs.
\newblock {\em Eur. J. Combin.} {\bf 16}, 405--414 (1995)

\bibitem{Urlep} Urlep, M.:
\newblock Triple intersection numbers of $Q$-polynomial distance-regular graphs.
\newblock {\em Eur. J. Combin.} {\bf 33(6)}, 1246--1252 (2012)

\bibitem{Vidali} Vidali, J.:
\newblock Using symbolic computation to prove nonexistence of distance-regular graphs.
\newblock \href{https://arxiv.org/abs/1803.10797}{arXiv:1803.10797}

\bibitem{Wang} Wang, H.C.:
\newblock Two-point homogeneous spaces.
\newblock {\em Ann. Math.} {\bf 55(2)}, 177--191 (1952)

\bibitem{Wilson} Wilson, R.M.:
\newblock An existence theory for pairwise balanced designs, III. Proof of the existence conjectures.
\newblock {\em J. Comb. Theory A} {\bf 18}, 71--79 (1975)

\bibitem{Aida} Yang, Q., Abiad, A., Koolen, J.H.:
\newblock An application of Hoffman graphs for spectral characterizations of graphs
\newblock {\em Electron. J. Comb.} {\bf 24(1)} P12 (2017)
 
\end{thebibliography}
\end{document}